\documentclass[11pt,twoside]{article}
\usepackage{amssymb,amsmath}
\usepackage{verbatim}
\usepackage[francais]{babel}
\def \P{{\hbox{\vrule width 0.6pt height 6.8pt depth -.2pt\kern-0.2pt
P}}}

\def \R {\mathbb R}
\def \T {\mathbb T}
\def \N {\mathbb N}
\def \e {\varepsilon}
\def\tilde{\widetilde}

\newcommand\Z{{\mathbb{Z}}}
\def\dq{\Delta_q}
\def\ddq{\dot\Delta_q}
\def\d{\partial}
\def\ov{\overline}

\def\cP{{\cal P}}
\def\cL{{\cal L}}
\def\cC{{\cal C}}
\def\cS{{\cal S}}

\def\div{\, \mbox{div}\,  }

\newenvironment{dem}{
\begin{description}
\item\textit{\textbf{D\'emonstration~:}}~}
{\hfill\rule{2.1mm}{2.1mm}
\end{description}}
\newcommand{\Sum}{\displaystyle \sum}

\newcommand{\Int}{\displaystyle \int}

\newcommand{\Sup}{\displaystyle \sup}

\newcommand{\du}{\delta\! u}
\newcommand{\dU}{\delta\! U}
\newcommand{\dt}{\delta\!\theta}
\newcommand{\dTheta}{\delta\!\Theta}
\newcommand\dPi{\delta\!\Pi}

\def \epsilon {\varepsilon}
\newenvironment{p}{
\begin{description}
\item\textit{\textbf{Preuve~:}}~}
{\hfill\rule{2.1mm}{2.1mm}
\end{description}}
\newtheorem{theo}{\bf TH\'EOR\`EME}[section]
\newtheorem{lem}[theo]{\bf LEMME}
\newtheorem{pro}[theo]{\bf PROPOSITION}

\newtheorem{defi}[theo]{\bf D\'EFINITION}
\newtheorem{rem}[theo]{\bf REMARQUE}

\begin{document}

\title{Les th\'eor\`emes de Leray
et de Fujita-Kato
 pour le syst\`eme de Boussinesq
partiellement visqueux.
\\
The Leray and Fujita-Kato theorems for the Boussinesq system
with partial viscosity}

\author{ R. Danchin et M. Paicu }
\maketitle
\centerline{\it Abstract}
\smallbreak
We are concerned with the so-called Boussinesq
equations with partial viscosity.
These equations consist of the ordinary incompressible
Navier-Stokes equations  with a forcing term
which is transported {\it with no dissipation} by the velocity field.
Such equations are simplified models for geophysics
(in which case the forcing term is  proportional
 either to the temperature, or to the salinity or to the density).

    In the present paper, we show that the standard theorems
 for incompressible Navier-Stokes equations may be extended to
 Boussinesq system despite the fact that there
is no dissipation or decay at large time for the forcing term.

More precisely, we state the global existence of finite
energy weak solutions in any dimension,
and global well-posedness in dimension $N\geq3$ for small data.
In the two-dimensional  case,
the finite energy global solutions are shown to be
unique for any  data in $L^2(\R^2).$

\centerline{\it R\'esum\'e} \smallbreak
Dans cet article, on
\'etudie le syst\`eme de Boussinesq  d\'ecrivant le ph\'enom\`ene de
convection  dans un fluide incompressible et visqueux. Ce syst\`eme
 est compos\'e des \'equations de
Navier-Stokes incompressibles avec un terme
de force verticale  dont l'amplitude est
transport\'ee {\it sans dissipation}
par le flot du champ de vitesses.

On montre que les r\'esultats classiques pour le syst\`eme
de Navier-Stokes standard demeurent vrais pour le syst\`eme
de Boussinesq bien qu'il n'y ait pas d'amortissement sur le
terme de force.

Plus pr\'ecis\'ement, on \'etablit l'existence de solutions faibles
globales d'\'energie finie en n'importe quelle dimension
et l'existence de solutions fortes uniques globales
en dimension $N\geq3$ pour de petites donn\'ees initiales.
Dans le cas particulier de la dimension deux, les solutions
d'\'energie finie sont uniques pour n'importe quelle donn\'ee
initiale dans $L^2(\R^2).$


 
\smallskip\noindent Mots-clefs: {\it Syst\`eme de Boussinesq, solutions faibles,
estimations \`a pertes, r\'egularit\'e critique.}

\section*{Introduction}

Dans cet  article, on \'etudie l'\'evolution d'un fluide incompressible
visqueux soumis \`a une force verticale.
On suppose que l'amplitude $\theta$ de cette force est
 transport\'ee par le flot du champ de vitesses $u$ du fluide.
Le  {\it  syst\`eme de Boussinesq} r\'egissant l'\'evolution
de $(\theta,u)$ s'\'ecrit donc~:
\begin{equation}\label{eq:boussinesq}
\begin{cases}
\partial_t\theta+u\cdot\nabla \theta=0\\
\partial_t u+u\cdot\nabla u-\nu \Delta u+\nabla\Pi=\theta\, e_N\\
\div u=0.
\end{cases}
\end{equation}
Les inconnues $\theta,$ $u$ et $\Pi$ d\'ependent
du temps $t\geq0$ et de la variable d'espace $x.$
Pour s'abstraire des difficult\'es dues aux conditions
aux limites, on suppose
que le fluide remplit tout l'espace (donc $x$ d\'ecrit $\R^N$ tout
entier).

Le syst\`eme \eqref{eq:boussinesq}
est un mod\`ele simplifi\'e courant pour
l'\'evolution de fluides g\'eophysiques (voir par exemple \cite{Ped}
ou \cite{Salmon}).
Le champ de vecteurs $u$ et le scalaire $\Pi$
 d\'esignent  alors respectivement  la vitesse et la pression
du fluide consid\'er\'e  et
 $\theta,$ une quantit\'e
scalaire\footnote{que nous appellerons temp\'erature pour fixer les
id\'ees mais qui  physiquement pourrait parfaitement
\^etre proportionnelle \`a
 la   densit\'e ou \`a  la salinit\'e} transport\'ee par le fluide.

Pour avoir un mod\`ele plus r\'ealiste, il conviendrait
de rajouter des termes de forces  ext\'erieures
aux deux \'equations.  Nous les supposerons nuls
pour simplifier la pr\'esentation.
Notons au passage que dans le cas particulier
o\`u $\theta\equiv0,$ on retrouve le syst\`eme de Navier-Stokes
``classique''.

Les pr\'emices de l'\'etude math\'ematique
du  syst\`eme de Boussinesq
sont relativement r\'ecentes et ne traduisent pas
uniquement une pr\'eoccupation d'ordre physique.
En r\'ealit\'e, l'engouement des math\'ematiciens
pour ce mod\`ele se limite essentiellement
au cas de la dimension deux du fait
d'une ressemblance formelle (dans le cas $\nu=0$)
avec  le syst\`eme d'Euler incompressible axisym\'etrique
{\it avec swirl}.
On peut  montrer que l'apparition de singularit\'es au temps $T$ est
li\'ee \`a l'explosion simultan\'ee de
$\nabla\theta$ et  du tourbillon dans $L^1(0,T;L^\infty)$
(voir \cite{ES}).
Malheureusement, d\'eterminer si ces quantit\'es explosent
effectivement semble au moins aussi difficile
que r\'epondre au  probl\`eme similaire pour
le syst\`eme d'Euler incompressible en dimension $3.$

Dans le cas $\nu>0$
en revanche,  les r\'esultats obtenus
sont nettement plus complets.
Divers auteurs ont \'etabli l'existence globale dans le cas
$N=2$ lorsque  l'\'equation sur $\theta$ comporte
en plus un terme de diffusion
(voir \cite{ST} et les r\'ef\'erences qui s'y trouvent).

L'existence globale en dimension deux {\it sans condition de petitesse}
demeure valable
dans la situation qui nous int\'eresse  o\`u $\theta$
est transport\'e sans diffusion.
Le cas de donn\'ees r\'eguli\`eres dans
des espaces de Sobolev a \'et\'e r\'esolu
(ind\'ependamment) par T. Hou et C. Li dans \cite{HL},
et par D. Chae dans \cite{CHAE}.
Tr\`es r\'ecemment, H. Abidi et T. Hmidi
ont \'etabli dans \cite{AH}
l'existence globale et unicit\'e sans condition de petitesse
pour des donn\'ees initiales $(\theta_0,u_0)$
appartenant \`a
 un tr\`es gros sous-espace de $(L^2(\R^2))^3.$
Comme on peut par ailleurs
construire des solutions faibles globales {\it \`a la Leray}
pour des donn\'ees initiales $\theta_0$
et $u_0$ dans  $L^2(\R^2)$
(voir par exemple \cite{HK1}), il est  raisonnable de
penser que
ces solutions ``faibles'' sont  uniques.
\smallbreak
Dans cet article, nous souhaitons  montrer que  trois des
r\'esultats math\'ematiques les plus c\'el\`ebres pour le
syst\`eme  de Navier-Stokes incompressible
sont encore valables pour le syst\`eme de
Boussinesq \eqref{eq:boussinesq}, \`a savoir~:
\begin{itemize}
\item l'existence de solutions faibles globales
d'\'energie finie en n'importe quelle dimension
(th\'eor\`eme de Leray, voir \cite{Leray}),
\item l'unicit\'e des solutions d'\'energie
finie en dimension deux (r\'esultat figurant implicitement dans
un autre article de J. Leray, voir \cite{Leray1},
et d\'emontr\'e par J.-L. Lions et G. Prodi dans \cite{LP}),
\item l'existence de solutions globales
uniques  pour des donn\'ees petites en dimension
$N\geq3$ (th\'eor\`eme de Fujita-Kato \cite{FK}).
\end{itemize}

Nous verrons plus loin
que la d\'emonstration  de l'existence globale de solutions faibles
se fait en suivant la d\'emarche originale de J. Leray~: la pr\'esence
 de la temp\'erature n\'ecessite un peu d'attention
 mais n'introduit pas
de difficult\'e majeure.

En revanche, pour d\'emontrer l'unicit\'e en dimension deux,
l'approche na\"ive consistant \`a
\'ecrire le syst\`eme v\'erifi\'e par la
diff\'erence de deux solutions semble vou\'ee \`a
l'\'echec. En effet,
du fait de la pr\'esence d'une \'equation
de transport dans le syst\`eme,
les estimations de stabilit\'e
ne peuvent \^etre \'etablies
qu'avec perte d'au moins une d\'eriv\'ee
donc dans un espace \`a r\'egularit\'e {\it strictement n\'egative}.
Or, pour contr\^oler ce type
de r\'egularit\'e, l'hypoth\`ese que
le champ de vitesses soit
au moins lipschitzien semble  in\'evitable.
Mais le champ de vitesses construit
  n'est a priori que  dans $L^2_{loc}(\R_+;H^1).$

Enfin, pour \'etablir l'existence de solutions fortes globales
en  dimension $N\geq3,$ une application directe de l'approche
 de Fujita-Kato associant
 point fixe contractant et
utilisation du semi-groupe de la chaleur ne saurait convenir.
 En effet, la m\'ethode de Fujita-Kato
 exige  non seulement  une donn\'ee initiale
petite mais aussi un terme de force $f$ qui soit
int\'egrable en temps sur $\R_+$ entier,
 et d'int\'egrale petite.
\`A moins que  $\theta$ ne soit identiquement nulle,
cette derni\`ere condition n'est  pas
v\'erifi\'ee par le terme de force $\theta\,e_N$
puisque $\theta$ est transport\'ee  sans amortissement
donc  est constante
le long des caract\'eristiques.
\medbreak
Nous avons adopt\'e le plan suivant.
La premi\`ere section est consacr\'ee \`a la pr\'esentation de
 nos r\'esultats principaux.
 Dans la partie suivante, nous d\'emontrons
l'analogue du th\'eor\`eme de Leray
pour le syst\`eme de Boussinesq
 \`a l'aide d'arguments \'el\'ementaires
d'analyse fonctionnelle.
Dans la troisi\`eme section,
nous pr\'esentons
l'artillerie technique n\'eces\-saire
pour venir \`a bout du reste de  notre programme~:
 espaces de Lorentz et de Besov, d\'ecom\-position
de Littlewood-Paley et calcul paradiff\'erentiel.
La quatri\`eme partie regroupe
diverses estimations a priori -- dont certaines
avec perte de r\'egularit\'e -- pour l'\'equation
de la chaleur ou de Stokes, l'\'equation de transport
et le syst\`eme de Navier-Stokes.
Les r\'esultats li\'es aux solutions fortes
sont alors d\'emontr\'es dans la cinqui\`eme partie, et l'unicit\'e
des solutions faibles d'\'energie finie
en dimension deux, dans la sixi\`eme partie.
Nous avons regroup\'e en derni\`ere section
 quelques  r\'esultats suppl\'ementaires.
Un lemme technique figure en appendice.


\section{R\'esultats principaux}\label{s:resultats}

Avant d'\'enoncer notre  r\'esultat d'existence globale
de solutions faibles, pr\'esentons bri\`evement
les \'egalit\'es d'\'energie formelles associ\'ees au syst\`eme.

Soit  donc $(\theta,u,\nabla\Pi)$ une  solution  r\'eguli\`ere
et d\'ecroissante \`a l'infini du syst\`eme
de Boussinesq avec donn\'ees $(\theta_0,u_0).$
Tout d'abord, sachant que le champ de vitesses est \`a divergence
nulle, on a
\begin{equation}\label{eq:temperature}
\Vert\theta(t)\Vert_{L^p}=\Vert\theta_0\Vert_{L^p}
\quad\text{pour tout}\ \ t\geq0\ \text{ et }\ 1\leq p\leq\infty.
\end{equation}
Par ailleurs, en prenant le produit scalaire $L^2(\R^N)^N$ de l'\'equation
de la vitesse avec $u$ puis en int\'egrant en temps,  on obtient
\begin{equation}\label{eq:vitesse}
\|u(t)\|^2_{L^2}+2\nu\int_0^t\|\nabla u\|^2_{L^2}\,d\tau=
\|u_0\|_{L^2}^2+2\int_0^t\int_{\R^N}\theta\,u_N\,dx\,d\tau.
\end{equation}
Il est clair que l'int\'egrale du  terme de droite peut \^etre absorb\'ee
par le membre de gauche
pourvu que l'on dispose d'un contr\^ole sur la norme $H^{-1}$
de $\theta.$ Ce contr\^ole est donn\'e
par \eqref{eq:temperature} d\`es que $L^p(\R^N)$ s'injecte contin\^ument
dans $H^{-1}(\R^N).$ Si $N\geq3$ (resp. $N=2$), cette condition
est v\'erifi\'ee si et seulement si  $p\geq\frac{2N}{N+2}$
(resp. $p>1$).

Comme de coutume, l'obtention des solutions faibles globales
se fera en passant \`a la limite dans une suite
de solutions approch\'ees v\'erifiant
\eqref{eq:temperature} et \eqref{eq:vitesse}.
 Le passage \`a la limite dans
le terme $u\cdot\nabla\theta=\div(u\theta)$
exige  que l'injection de  $L^p(\R^N)$
dans $H^{-1}(\R^N)$ soit {\it localement compacte},
donc que  $p>\frac{2N}{N+2}.$

Ces consid\'erations sommaires m\`enent \`a l'\'enonc\'e suivant~:
\begin{theo}\label{solutions-faibles-energie-finie}
Soit $\theta_0\in L^p(\R^N)$ avec $\frac{2N}{N+2}< p\leq
 2$,
 et $u_0$ un champ de vitesses \`a divergence nulle
et coefficients dans $L^2(\R^N).$
 Alors le syst\`eme de Boussinesq avec donn\'ee initiale
$(\theta_0,u_0)$ admet une solution
faible globale $$(\theta, u)\in L^\infty(\R_+;L^p(\R^N))\times
\big(L^\infty_{loc}(\R_+;L^2(\R^N))\cap L^2_{loc}(\R_+;H^1(\R^N))\big)^N,$$
telle que  pour tout
$t\geq0,$ on ait
$$\|\theta (t)\|_{L^p}\leq
\|\theta_0\|_{L^p}\quad\text{et}\quad
\|u(t)\|^2_{L^2}+\nu\int_0^t\|\nabla u\|^2_{L^2}d\tau \leq
C\Bigl(\|u_0\|_{L^2}^2+\nu^{\alpha-1}t^{\alpha+1}\|\theta_0\|_{L^p}^2
\Bigr)$$
pour une constante $C$ ne d\'ependant que de $N$ et de $p,$
et $\alpha=1-N\bigl(\frac1p-\frac12\bigr).$
\end{theo}
Dans le cas $\theta\equiv0,$
il est bien connu que les solutions faibles
globales en dimension deux sont uniques (voir par
exemple \cite{LP}).
Il est donc l\'egitime de chercher \`a d\'emontrer
que cette propri\'et\'e persiste
pour le syst\`eme de Boussinesq avec des donn\'ees
quelconques dans $L^2(\R^2).$
\`A notre connaissance, le r\'esultat le plus proche
allant dans ce sens a \'et\'e \'etabli
par H. Abidi et T. Hmidi dans \cite{AH}
(voir aussi \cite{HK1})~: si
$\theta_0$ appartient \`a l'espace
de Besov homog\`ene\footnote{voir la d\'efinition \ref{def}}
 $\dot B^0_{2,1}(\R^2)$ et
 $u_0$ est dans $\dot B^{-1}_{\infty,1}\cap
L^2(\R^2)$ alors  il y a unicit\'e de la solution globale
construite.

Nous d\'emontrons ici qu'il y a encore unicit\'e
des solutions faibles d'\'energie finie en dimension
deux {\it sous la seule hypoth\`ese que les
donn\'ees sont dans $L^2(\R^2).$}

Cette propri\'et\'e remarquable
d\'ecoule du fait que
le champ de vitesses construit appartient aussi \`a
l'espace  $\tilde L^1_{loc}(\R_+;H^2(\R^2))$
(d\'efini plus loin en \eqref{eq:ltilde})
qui est \`a peine plus gros que
$L^1_{loc}(\R_+;H^2(\R^2)).$
De ce fait, le champ de vitesses
est presque localement int\'egrable en temps
\`a valeurs quasi-lipschitziennes et l'on peut faire appel
\`a des estimations avec perte de r\'egularit\'e
dans l'esprit de celles qui ont \'et\'e
\'etablies dans \cite{BCh} ou, plus r\'ecemment, dans \cite{DANCHI}.
\begin{theo}\label{th:trois}
 Supposons $N=2.$ Soit $(\theta_0,u_0)\in L^2(\R^2)^3$ avec
$u_0$ \`a divergence nulle.
 Alors le syst\`eme $\eqref{eq:boussinesq}$ a une unique
solution globale $(\theta,u,\nabla\Pi)$ telle que
$$
\theta\in\cC(\R_+;L^2)\quad\text{et}\quad
u\in\cC(\R_+;L^2)\cap L^2_{loc}(\R_+;H^1).
$$
De plus la norme $L^2$ de $\theta$ est conserv\'ee
au cours de l'\'evolution, l'\'egalit\'e
 $\eqref{eq:vitesse}$ est satisfaite et $u$
appartient \`a l'espace $\tilde L^1_{loc}(\R_+;H^2)$
d\'efini au-dessus de la remarque $\ref{r:stokes}$.
\end{theo}

La derni\`ere partie de notre programme consiste \`a \'etablir
un r\'esultat d'existence globale de solutions fortes
\`a donn\'ees petites en dimension $N\geq3,$
dans l'esprit de celui de  Fujita-Kato pour le
syst\`eme de Navier-Stokes incompressible (voir \cite{FK}).

Rappelons que le th\'eor\`eme de Fujita-Kato
se d\'emontre en r\'ecrivant  le syst\`eme
de Navier-Stokes en termes de probl\`eme de  point fixe
pour une fonctionnelle  construite
\`a l'aide du semi-groupe de Stokes.
Sous des hypoth\`eses de petitesse ad\'equates, le
th\'eor\`eme du point fixe de Picard permet alors
d'obtenir  une  solution globale unique.
Cette approche  s'av\`ere particuli\`erement performante
si l'espace fonctionnel $F$ utilis\'e
est {\it critique}, i.e.
 respecte l'invariance par changement d'\'echelle
(ou {\it scaling} en anglais)  du syst\`eme de Navier-Stokes.
Cela am\`ene \`a choisir pour $F$ un espace de distributions
sur $\R^+\times\R^N$ \`a norme invariante pour tout
$\lambda>0$ par la transformation
$$
u(t,x)\longmapsto\lambda u(\lambda^2t,\lambda x)
$$
et donc \`a choisir une vitesse initiale $u_0$
dans un espace fonctionnel $E$ \`a norme invariante
par la transformation  $u_0\mapsto\lambda u_0(\lambda\,\cdot\,).$
 H. Fujita et T. Kato d\'emontrent
ainsi que le syst\`eme de Navier-Stokes est globalement
bien pos\'e pour des donn\'ees initiales petites
par rapport \`a la viscosit\'e dans  l'espace 
de Sobolev homog\`ene $\dot H^{\frac12}(\R^3)$
(voir \cite{FK} et \cite{Che}).  
Le r\'esultat de
Fujita-Kato se g\'en\'eralise  \`a un grand nombre
d'espaces critiques
(voir e.g. \cite{L} ou \cite{Meyer}).
\smallbreak
 Un calcul facile montre que  le syst\`eme de Boussinesq est invariant
par  la transformation
$$
u(t,x)\longmapsto\lambda u(\lambda^2t,\lambda x)\quad\text{et}\quad
\theta(t,x)\longmapsto\lambda^3\theta(\lambda^2t,\lambda x).
$$
Autrement dit, les espaces critiques
pour la vitesse sont les m\^emes que pour le syst\`eme de
Navier-Stokes et il faut  en quelque sorte exiger
{\it deux d\'eriv\'ees de moins} sur la temp\'erature.
\smallbreak
En l'absence de dissipation sur la temp\'erature,
il ne semble pas possible de
d\'emontrer que le syst\`eme de Boussinesq est bien pos\'e
m\^eme localement pour des donn\'ees aussi peu r\'eguli\`eres.
On peut cependant \'etablir
l'existence locale pour des donn\'ees
initiales dans l'espace de Besov homog\`ene $\dot B^0_{N,1}(\R^N)$
(gros sous-espace de $L^N(\R^N)$)
qui est critique pour la vitesse seulement~:
\begin{theo}\label{th:quatre} Supposons $N\geq2.$
Soit $\theta_0\in \dot B^0_{N,1}$
et $u_0$ un champ de vecteurs \`a divergence nulle et coefficients
dans $\dot B^0_{N,1}.$ Alors le syst\`eme de Boussinesq admet une unique
solution locale $(\theta,u,\nabla\Pi)$ dans l'espace
$$E_T:=\cC([0,T];\dot B^0_{N,1})\times
\Bigl(\cC([0,T];\dot B^0_{N,1})\cap L^1(0,T;\dot B^2_{N,1})\Bigr)^N
\times\Bigl(L^1(0,T;\dot B^0_{N,1})\Bigr)^N.
$$
Si de plus $\,\theta\in L^\infty(0,T^*;\dot B^0_{N,1})\,$
et $\,u\in L^\infty(0,T^*;\dot B^0_{N,1})\cap L^1(0,T^*;\dot B^2_{N,1})\,$
alors la solution peut \^etre prolong\'ee au-del\`a de $T^*.$
\end{theo}
Si la construction de solutions \`a temps petit peut
se faire par des arguments standard,
montrer l'existence de solutions fortes globales
est nettement plus d\'elicat.
En effet, comme expliqu\'e dans
l'introduction,  l'absence  d'amortissement pour la temp\'erature
interdit l'approche classique associant
estimations pour le  semi-groupe de la chaleur et point
fixe de Picard, et le scaling de la temp\'erature semble
bien trop bas pour pouvoir travailler directement dans
des espaces invariants par changement d'\'echelle.

Notre strat\'egie  est la suivante~: trouver deux espaces
fonctionnels, $X$ pour la vitesse initiale $u_0$ et
$Y$ pour la temp\'erature initiale $\theta_0$
 tels que, en notant $\cP$ le projecteur
orthogonal sur les champs \`a divergence nulle et
\ $\bigl(e^{\tau\Delta}\bigr)_{\tau>0}$ le
semi-groupe de la chaleur, on ait
 les estimations suivantes:
$$\begin{array}{c}\|\theta\|_{L^\infty(0,T;Y)}\leq \|\theta_0\|_{Y},
\qquad\|e^{t\Delta}u_0\|_{L^\infty(0,T;X)}\leq \|u_0\|_{X},\\[1.5ex]
\bigg\|\Int_0^t e^{(t-s)\Delta}\cP(\theta
e_N)\,ds\bigg\|_{L^\infty(0,T;X)}\leq  C_1\|\theta
\|_{L^\infty(0,T;Y)},\\[2.5ex] \bigg\|\Int_0^t e^{(t-s)\Delta}
\cP\div(u\otimes u)(s)\,ds\bigg\|_{L^\infty(0,T;X)} \leq
C_2\|u\|^2_{L^\infty(0,T;X)}\end{array}$$
avec des constantes $C_1$ et $C_2$
{\it ind\'ependantes} du temps $T.$ (Nous avons
suppos\'e que  $\nu=1$ pour simplifier l'heuristique.)
\smallbreak
Pour de tels espaces, on obtiendra alors
$$\|u\|_{L^\infty(0,T;X)}\leq
\|u_0\|_{X}+C_1\|\theta_0\|_{Y}+C_2\|u\|^2_{L^\infty(0,T;X)}.$$
Supposons que
\eqref{eq:boussinesq} admette une solution locale
$(\theta,u)$ dans $L^\infty([0,T];Y)\times L^\infty([0,T];X)$
et que de plus
  $C_1\|\theta_0\|_{Y}+\|u_0\|_{X}\leq c$ avec $c$ suffisamment petit.
Soit $T^*$  le temps maximal  pour lequel
$\|u\|_{L^\infty(0,T^*;X)}\leq 2c.$
Si l'on a choisi $c$ de telle sorte que $4cC_2<1,$
l'in\'egalit\'e  ci-dessus implique que $\|u\|_{L^\infty(0,T^*;X)}<2c.$
Par des arguments classiques de bootstrap,
  on conclut alors que $T^*=+\infty.$

Nous verrons plus loin que l'on peut choisir pour  $X$ l'espace
de Lorentz $L^{N,\infty}(\R^N),$
et pour $Y$ l'espace de Lebesque $L^{\frac N3}(\R^N)$
(voire  $L^{\frac N3,\infty}(\R^N)$ si $N\geq4$).
Notons au passage que ces espaces sont {\it critiques}.

\`A moins que $\theta\equiv0,$
il n'est pas clair que l'on puisse
d\'emontrer l'existence et l'unicit\'e de  solutions
dans ces espaces sans hypoth\`ese
suppl\'ementaire de r\'egularit\'e.
Pour y rem\'edier, nous allons devoir  exiger
plus de r\'egularit\'e sur les donn\'ees initiales.
Cette r\'egularit\'e doit de plus
pouvoir \^etre transport\'ee
{\it globalement} si les donn\'ees
sont petites dans $X\times Y.$
Le caract\`ere global demande un peu
d'attention car les estimations habituelles de r\'egularit\'e
sur les solutions d'une \'equation de transport
 font appara\^{\i}tre un  terme
qui d\'epend {\it exponentiellement} de la norme lipschitzienne
du champ de vitesse.
On sait cependant depuis un article
de M. Vishik que  ces estimations
sont ``meilleures'' dans les espaces de
Besov d'indice nul  (voir \cite{Vishik} et la proposition
\ref{p:transport2} pour plus de d\'etails).
\smallbreak
Ces consid\'erations m\`enent au  r\'esultat suivant
d'existence globale \`a donn\'ees petites~:
\begin{theo}\label{th:deux} Supposons que $N\geq3.$
Soit $\theta_0\in \dot B^0_{N,1}(\R^N)$
 et $u_0\in \dot B^{-1+\frac Np}_{p,1}(\R^N)$
 pour un $p\in[N,\infty].$
Il existe une constante $c>0$ ne d\'ependant que de $N$
telle que si de plus $u_0\in L^{N,\infty},$ $\theta_0\in L^{\frac N3}$
et
\begin{equation}\label{eq:petit}
\|u_0\|_{L^{N,\infty}}+\nu^{-1}\|\theta_0\|_{L^{\frac N3}}\leq c\nu
\end{equation}
alors le syst\`eme de Boussinesq admet une unique solution globale
$$
(\theta,u,\nabla\Pi)\in{\cal C}(\R_+;\dot B^0_{N,1})\times
\Bigl({\cal C}(\R_+;\dot B^{\frac Np-1}_{p,1}) \cap
L^1_{loc}(\R_+;\dot B^{\frac Np+1}_{p,1})\Bigr)^N \times
\Bigl(L^1_{loc}(\R_+;\dot B^{\frac Np-1}_{p,1})\Bigr)^N.
$$
De plus, il existe une constante $C$ ne d\'ependant
que de $N$ telle que pour tout temps $t\geq0,$ on~ait
$$
\displaylines{\|\theta (t)\|_{L^{\frac N3}}=\|\theta_0\|_{L^{\frac
N3}}, \qquad \|u(t)\|_{L^{N,\infty}}\leq C\Bigl(
\|u_0\|_{L^{N,\infty}}+\nu^{-1}\|\theta_0\|_{L^{\frac N3}}\Bigr),\cr
\|u\|_{L^\infty_t(\dot B^{\frac Np-1}_{p,1})}
+\nu\|u\|_{L^1_t(\dot B^{\frac Np+1}_{p,1})}\leq
C\|u_0\|_{\dot B^{\frac Np-1}_{p,1}}e^{Ct\nu^{-1}\|\theta_0\|_{\dot B^0_{N,1}}}
+\nu\Bigl(e^{Ct\nu^{-1}\|\theta_0\|_{\dot B^0_{N,1}}}-1\Bigr),\cr
\|\theta(t)\|_{\dot B^0_{N,1}}\leq
C\|\theta_0\|_{\dot B^0_{N,1}}e^{Ct\nu^{-1}\|\theta_0\|_{\dot B^0_{N,1}}}
\bigl(1+C\nu^{-1}\|u_0\|_{\dot B^0_{N,1}}\bigr).}$$
\end{theo}
\begin{rem} \`A partir de la dimension quatre,
l'espace $L^{\frac N3}(\R^N)$ peut \^etre remplac\'e par $L^{\frac
N3,\infty}(\R^N).$
\end{rem}
\begin{rem}
Notons que les hypoth\`eses suppl\'ementaires
sur la vitesse sont critiques en terme de scaling.
Par ailleurs,  l'appartenance de $\theta_0$
\`a l'espace de Besov homog\`ene $\dot B^0_{N,1}$
(qui est un sous-espace
strict de l'espace de Besov non homog\`ene $B^0_{N,1}$) est
assur\'ee par $\theta_0\in L^{\frac N3,\infty}\cap B^0_{N,1}.$
De m\^eme, si $p>N,$ l'hypoth\`ese $u_0\in \dot B^{-1+\frac Np}_{p,1}$
est garantie par $u_0\in L^{N,\infty}\cap B^{-1+\frac
Np}_{p,1}$ (voir le lemme $\ref{l:lorentz}$).
\end{rem}
\noindent {\bf Notation~:}
Dans tout l'article, $C$ d\'esigne une ``constante''
susceptible de changer de ligne en ligne
et dont la valeur exacte n'influe pas sur l'exactitude des
calculs. On utilise parfois la  notation abr\'eg\'ee
$A\lesssim B$ au lieu de $A\leq CB.$


\section{Solutions faibles globales}\label{s:faible}

Cette section est consacr\'ee \`a la d\'emonstration
du th\'eor\`eme \ref{solutions-faibles-energie-finie}.

Fixons une fonction $\chi\in\cC_c^\infty(\R^N)$
positive et d'int\'egrale $1,$
et notons $I_r$ ($r>0$) l'op\'erateur de convolution
par  la fonction $r^{-N}\chi(r^{-1}\cdot).$
Dans un premier temps, nous allons 
r\'esoudre le syst\`eme de Boussinesq  r\'egularis\'e 
suivant~: 
\begin{equation}\label{eq:regularise}
\begin{cases}
\partial_t\theta+I_ru\cdot\nabla\theta=0,\\
\partial_tu+\cP\bigl(I_ru\cdot\nabla u\bigr) 
-\nu\Delta u={\cal P}\bigl(I_r\theta\, e_N\bigr),
\end{cases}
\end{equation}
avec donn\'ees initiales $\theta_0$ et $u_0$ dans $L^2(\R^N).$
\begin{pro} 
Pour tout couple de donn\'ees initiales 
$(\theta_0,u_0)$ dans $L^2(\R^N)$ tel que ${\rm div}\, u_0=0,$
le syst\`eme $\eqref{eq:regularise}$ admet 
une solution faible globale $(\theta,u)$
telle que 
$$
\theta\in\cC(\R_+;L^2),\qquad
u\in L^\infty_{loc}(\R_+;L^2)\cap L^2_{loc}(\R_+;H^1)
$$
et v\'erifiant l'in\'egalit\'e d'\'energie suivante pour
tout $t\in\R^+$:
\begin{equation}\label{eq:vitesse1}
\|u(t)\|^2_{L^2}+2\nu\int_0^t\|\nabla u\|^2_{L^2}\,d\tau\leq
\|u_0\|_{L^2}^2+2\int_0^t\int_{\R^N}I_r\theta\,u_N\,dx\,d\tau.
\end{equation}
\end{pro}
\begin{p}
Nous allons r\'esoudre
le syst\`eme \eqref{eq:regularise}
\`a l'aide de
  la m\'ethode de Friedrichs. Pour cela,
on d\'efinit l'op\'erateur  de troncature spectrale $J_n$ par
$$\widehat{J_nf}(\xi)=1_{[\frac 1n,n]}(|\xi|)\widehat
f(\xi)$$
et l'on cherche \`a r\'esoudre
 le syst\`eme  suivant pour $n\geq1$:
\begin{equation}\label{approche}\begin{cases}
\partial_t\theta^n+J_n\bigl(I_r{\cal P}J_nu^n\cdot\nabla J_n\theta^n\bigr)=0,\\
\partial_tu^n+\cP J_n\bigl(I_r{\cal P}J_nu^n\cdot\nabla{\cal P}J_nu^n\bigr) 
-\nu\Delta{\cal P}J_nu^n={\cal P}\bigl(I_rJ_n\theta^n\, e_N\bigr),\\
(\theta^n,u^n)_{|t=0}=(J_n\theta_0,J_n u_0).
\end{cases}
\end{equation}
Ce syst\`eme est une \'equation
diff\'erentielle ordinaire sur $L^2$ v\'erifiant les hypoth\`eses
du th\'eor\`eme de Cauchy-Lipschitz.
Il admet donc une unique solution maximale 
$(\theta^n,u^n)\in\cC^1\bigl([0,T_n[;(L^2(\R^N))^{N+1}\bigr)$
pour un $T_n>0.$

  Comme de plus $(\cP J_n)^2=\cP J_n$ et $J_n^2=J_n,$
le couple $(J_n\theta^n,\cP J_nu^n)$ est aussi solution 
de \eqref{approche}. Par unicit\'e, on~a donc
$J_n\theta^n=\theta^n$ et  $\cP J_nu^n=u^n.$
 Cela entra\^{\i}ne
 que  $\theta^n$ et $u^n$ sont dans $\cC^1([0,T^n[;H^\infty)$
et v\'erifient
\begin{equation}\label{eq:approche1}
\begin{cases}
\partial_t\theta^n+J_n\bigl(I_ru^n\cdot\nabla\theta^n\bigr)=0,\\
\partial_tu^n+\cP J_n\bigl(I_ru^n\cdot\nabla u^n\bigr) 
-\nu\Delta u^n={\cal P}\bigl(I_r\theta^n\, e_N\bigr).\end{cases}
\end{equation}
Une m\'ethode d'\'energie \'el\'ementaire assure alors que
 $\|\theta^n(t)\|_{L^2}=\|J_n\theta_0\|_{L^2}$ pour
tout $t\in\R_+,$ et
\begin{equation}\label{eq:vitesse2}
\|u^n(t)\|^2_{L^2}+2\nu\int_0^t\|\nabla u^n\|^2_{L^2}\,d\tau=
\|J_nu_0\|_{L^2}^2+2\int_0^t\int_{\R^N}I_r\theta^n\,u^n_N\,dx\,d\tau.
\end{equation}
Notons que ces deux \'egalit\'es
impliquent que $T_n=+\infty.$
En effet, si $T_n$ est fini, 
on montre facilement en majorant le terme
de droite de \eqref{eq:vitesse2}
\`a l'aide de l'in\'egalit\'e de Cauchy-Schwarz, 
que  $(\theta^n,u^n)$ est dans  $L^\infty(0,T^n;L^2)$
et  peut donc \^etre prolong\'ee
en vertu des th\'eor\`emes classiques sur les \'equations
diff\'erentielles ordinaires.
\smallbreak
Il s'agit maintenant de passer \`a la limite 
dans \eqref{eq:approche1}.
Tout d'abord, il est clair que la suite $(\theta^n)_{n\in\N}$
est born\'ee dans $L^\infty(\R_+;L^2)$
et que \eqref{eq:vitesse2} assure que 
$(u^n)_{n\in\N}$ est born\'ee dans $L^\infty_{loc}(\R_+;L^2)
\cap L^2_{loc}(\R_+;H^1).$
En utilisant le syst\`eme \eqref{eq:approche1} et en 
remarquant que la suite $(I_ru^n)_{n\in\N}$ est born\'ee dans tous les espaces
$L^\infty_{loc}(\R_+;H^s),$ 
on en d\'eduit alors que 
$(\d_t\theta^n)_{n\in\N}$ et $(\d_tu^n)_{n\in\N}$
sont born\'ees dans $L^2_{loc}(\R_+;H^{-1}).$
Par  des arguments  de convergence faible et de compacit\'e classiques,
  on conclut  qu'il 
existe $\theta\in L^\infty(\R_+;L^2)$
et $u\in L^\infty_{loc}(\R_+;L^2)\cap L^2_{loc}(\R_+;H^1)$
\`a divergence nulle, tels que, \`a extraction pr\`es,\smallbreak
\begin{itemize}
\item $\theta^n$ tend vers $\theta$
dans $L^\infty_{loc}\bigl(\R_+;H^{-\eta}_{loc}\bigr)$
 pour tout $\eta>0,$
\item $\theta^n(t)$ tend faiblement dans $L^2$ vers
$\theta(t)$ pour tout $t\in\R_+$ (et donc
$I_r\theta^n(t)$ tend {\it fortement}  vers
$I_r\theta(t)$ dans $L^2$),
\smallbreak
\item  $u^n$ tend vers  $u$  dans
$L^{2}_{loc}(\R_+;H^{1-\epsilon}_{loc})$ pour tout $\epsilon>0,$
\item  $u^n(t)$ tend faiblement dans $L^2$ vers
$u(t)$ pour tout $t\in\R_+.$
\end{itemize}\smallbreak
Ces propri\'et\'es de convergence  permettent 
de passer \`a la limite 
dans le syst\`eme \eqref{eq:approche1}
et de montrer que  
$$
\lim_{n\rightarrow\infty}
\int_0^t\int_{\R^N}I_r\theta^n\,u^n_N\,dx\,d\tau
=\int_0^t\int_{\R^N}I_r\theta\,u_N\,dx\,d\tau.
$$ 
Le couple $(\theta,u)$ est
donc une solution globale 
de \eqref{eq:regularise} v\'erifiant
\eqref{eq:vitesse1}.
\end{p} 
Nous pouvons maintenant passer \`a la 
 d\'emonstration
du th\'eor\`eme \ref{solutions-faibles-energie-finie}.
Nous supposons donc que la temp\'erature initiale
$\theta_0$ appartient \`a $L^p(\R^N)$ avec $\frac{2N}{N+2}<p\leq2$
et que $u_0\in L^2(\R^N)$ avec $\div u_0=0.$
Il est clair que $I_r\theta_0\in L^2$ pour tout $r>0.$
La proposition pr\'ec\'edente appliqu\'ee avec $r=2^{-k}$
assure donc l'existence d'une
solution globale $(\theta_k,u_k)$ 
v\'erifiant \eqref{eq:vitesse1}
pour le syst\`eme approch\'e
\begin{equation}\label{eq:approche2}
\begin{cases}
\partial_t\theta_k+\div\bigl(I_{2^{-k}}u_k\,\theta_k\bigr)=0,\\
\partial_tu_k+\cP\div\bigl(I_{2^{-k}}u_k\otimes u_k\bigr) 
-\nu\Delta u_k={\cal P}\bigl(I_{2^{-k}}\theta_k\, e_N\bigr),\\
(\theta_k,u_k)_{|t=0}=(I_{2^{-k}}\theta_0,u_0).
\end{cases}
\end{equation}
On constate que $\theta_k$
est solution d'une \'equation de transport 
par le champ de vecteurs  r\'egulier \`a divergence nulle
$I_{2^{-k}}u_k.$ Cela assure que 
$\theta_k\in\cC(\R_+;L^2\cap L^p)$ et que 
 \begin{equation}\label{eq:theta1}
\|\theta_k(t)\|_{L^p}= \|I_{2^{-k}}\theta_0\|_{L^p}\leq
 \|\theta_0\|_{L^p}
\quad\text{pour tout}\quad t\geq0.
\end{equation}
Ensuite, en vertu de \eqref{eq:vitesse1}
 et de l'in\'egalit\'e de H\"older, on~a
 \begin{equation}\label{energie}
\|u_k(t)\|^2_{L^2}+2\nu\int_0^t\|\nabla u_k\|^2_{L^2}\,dt
\leq \|u_0\|^2_{L^2}
+ 2\int_0^t\|\theta_k\|_{L^p}\|u_k\|_{L^{p'}}\,dt.\end{equation}
Comme $\frac{2N}{N+2}<p\leq 2$,  on dispose de l'in\'egalit\'e de
 Gagliardo-Nirenberg suivante~:
$$\|u_k\|_{L^{p'}}\leq
C\|u_k\|_{L^2}^{\alpha}\|\nabla u_k\|_{L^2}^{1-\alpha}
\quad\text{avec}\quad \alpha=1-N\Bigl(\frac{1}{p}-\frac{1}{2}\Bigr).
$$
En injectant cette in\'egalit\'e dans
\eqref{energie} et en
utilisant \eqref{eq:theta1} et  l'in\'egalit\'e de
Young, on trouve
$$
\|u_k(t)\|^2_{L^2}+\nu\int_0^t\|\nabla u_k\|^2_{L^2}\,dt
\leq \|u_0\|^2_{L^2}
+C\nu^{\frac{\alpha-1}{\alpha+1}}\|\theta_0\|_{L^p}^{\frac2{1+\alpha}}
\int_0^t\|u_k\|_{L^2}^{\frac{2\alpha}{1+\alpha}}\,d\tau,
$$
d'o\`u, apr\`es une nouvelle utilisation de l'in\'egalit\'e
de Young, 
\begin{equation}\label{energie-2}
\|u_k(t)\|^2_{L^2}+\nu\int_0^t\|\nabla u_k\|^2_{L^2}\,dt
\leq C\bigl(\|u_0\|^2_{L^2}
+t^{\alpha+1}\nu^{\alpha-1}\|\theta_0\|_{L^p}^2\bigr).
\end{equation}
Il  reste \`a justifier  le passage \`a la limite
dans \eqref{eq:approche2}.
Tout d'abord \eqref{eq:theta1} et \eqref{energie-2}
assurent que $(\theta_k)_{k\in\N}$ est 
born\'ee dans $L^\infty(\R_+;L^p)$
et que $(u_k)_{k\in\N}$ est
born\'ee dans $L^\infty(\R_+;L^2)\cap L^2(\R_+;H^1).$
En  utilisant \eqref{approche}, on peut
alors montrer que
$(\d_t\theta_k)_{k\in\N}$
et   $(\d_tu_k)_{k\in\N}$ sont  born\'ees dans
$L^2_{loc}(\R_+;H^{-\frac N2-1}).$
Par  des arguments  de compacit\'e classiques
et en notant que $L^p(\R^N)\hookrightarrow H^{-N(\frac
1p-\frac12)}(\R^N)$ pour $1<p\leq2,$
  on en d\'eduit, quitte \`a extraire, que~:
\begin{itemize}
\item $\theta_k$ tend vers une fonction $\theta$
dans $L^\infty_{loc}\bigl(\R_+;H^{-N(\frac 1p-\frac12)-\eta}_{loc}\bigr)$
 pour tout $\eta>0,$
\smallbreak
\item $u_k$ tend vers un champ $u$ \`a divergence nulle dans
$L^{2}_{loc}(\R_+;H^{1-\epsilon}_{loc})$ pour tout $\epsilon>0.$
\end{itemize}\smallbreak
Ces propri\'et\'es de convergence permettent
de justifier que $\div(I_{2^{-k}}u_k\theta_k)$
tend vers $\div(u\theta)$ au sens des distributions
pourvu que $1-\frac Np+\frac N2>0,$
 condition qui est \'equivalente \`a $p>2N/(N+2).$
Le passage \`a la limite dans le terme $\div(I_{2^{-k}}u_k\otimes u_k)$
se fait exactement comme pour le syst\`eme
de Navier-Stokes incompressible ``classique".


\section{Espaces fonctionnels et outils d'analyse harmonique}

Commen\c cons par rappeler la d\'efinition
des espaces de type $L^p$ faible.
\begin{defi} Pour $1\leq p<\infty,$ on d\'efinit
l'espace $L^p$ faible (not\'e $L^{p,\infty}$ par la suite)
 comme l'ensemble des fonctions
mesurables de $\R^N$ dans $\R$ telles que
$$
\|f\|_{L^{p,\infty}}:=\sup_{\lambda>0} \lambda\:
\Bigl|\bigl\{x\in\R^N/|f(x)|>\lambda\bigr\}\Bigr|^{\frac1p}<\infty.
$$
\end{defi}
\begin{rem}\label{r:lorentz}
Dans le cas $1<p<\infty,$ l'espace $L^{p,\infty}$
co\"{\i}ncide avec l'espace
de Lorentz d\'efini par l'interpolation r\'eelle
$(L^\infty,L^1)_{(\frac{1}{p},\infty)}.$ En d'autres termes,
 $f$ est dans $L^{p,\infty}$ si et seulement si, pour
tout $A>0$, on peut \'ecrire $f=f^A+f_A$ avec $\|f_A\|_{L^1}\leq
CA^{1-1/p}$ et $\|f^A\|_{L^\infty}\leq C A^{-1/p},$
et la ``meilleure constante" $C$ est  une norme
\' equivalente \`a la quantit\'e de la d\'efinition
 ci-dessus (voir par exemple \cite{L}).

Plus g\'en\'eralement, pour $1\leq q\leq\infty,$ l'espace de Lorentz
$L^{p,q}$ peut \^etre d\'efini par
interpolation r\'eelle en posant~:
$$L^{p,q}:=
(L^\infty,L^1)_{(\frac{1}{p},q)}.$$
\end{rem}
Dans la suite de cette partie, nous rappelons bri\`evement la
d\'efinition de la
d\'ecomposition de Littlewood-Paley (selon \cite{C}) ainsi que celle
des espaces de Besov.

La d\'ecomposition de  Littlewood-Paley
se d\'efinit \`a l'aide  d'une d\'ecomposition dyadique de l'unit\'e~:
soit  $\chi$ une fonction
positive radiale  de classe $\cC^\infty,$
support\'ee dans la boule $\{\vert\xi\vert\leq\frac{4}{3}\},$
valant $1$ sur  $\{\vert\xi\vert\leq\frac{3}{4}\},$
et telle que $r\mapsto \chi(re_r)$ soit d\'ecroissante.
On pose  $\varphi(\xi)=\chi(\xi/2)-\chi(\xi)$ de telle sorte que
\begin{equation}\label{eq:dyadique}
\forall\xi\in\R^N\setminus\{0\},\;
\sum_{q\in \Z}\varphi(2^{-q}\xi)=1\quad\text{et}\quad
\forall\xi\in\R^N,\; \chi(\xi)=1-\sum_{q\in\N}\varphi(2^{-q}\xi).
\end{equation}
On d\'efinit ensuite les op\'erateurs de localisation en fr\'equences
$\dot\Delta_q$ et $\dot S_q$ de ${\cal{L}}({\cal{S}}';\cS)$ par:
$$
\dot\Delta_q\;u:=\varphi(2^{-q}D)u\quad\text{et}\quad
\dot S_q\,u:=\chi(2^{-q}D)u\quad \mbox{pour tout}\quad q\in\Z.
$$
 Notons que pour une distribution temp\'er\'ee
$u\in\mathcal S'(\R^N)$, les  fonctions $\dot\Delta_q u$ et
$\dot S_qu$ sont analytiques \`a
croissance lente. Si de plus il existe un r\'eel $s$ tel que $u\in
H^s(\R^N)$ alors  $\dot\Delta_q u$ et $\dot S_qu$  appartiennent \`a
l'espace
$H^\infty:=\bigcap_{\sigma\in\R} H^\sigma.$

 Par ailleurs, on
\'etablit  \`a l'aide de \eqref{eq:dyadique}
 que
$\ \displaystyle{
u=\dot S_0u+\sum_{q\in\N}\dot\Delta_qu}\ $
dans $\cS'(\R^N)$ et
$\ \displaystyle{u=\sum_{q\in\Z}\dot\Delta_q u}\ $ dans
${\cal {S}}'(\R^N)$ modulo un polyn\^ome (voir par exemple \cite{L}).
\medbreak
 Enfin, nous utiliserons souvent  la propri\'et\'e de presque
orthogonalit\'e suivante~:
\begin{equation}\label{presortho}
\dot\Delta_k\dot\Delta_q u\equiv 0 \quad\mbox{si}\quad\vert k-q\vert\geq 2
\quad\mbox{et}\quad\dot\Delta_k( \dot S_{q-1}u\dot\Delta_qv)\equiv
0\quad\mbox{si} \quad\vert k-q\vert\geq 5.
\end{equation}
Nous pouvons maintenant d\'efinir les espaces de Besov homog\`enes
comme suit~:
\begin{defi}\label{def}
Soient $s\in\R,\,(p,r)\in[1,\infty]^2$ et $u\in {\cal {S}}'(\R^N).$
On note
$$
\Vert u\Vert_{ \dot B^s_{p,r}} \overset{\rm d\acute{e}f}{=}
\Big(\sum_{q\in\Z}2^{rqs}
\Vert\dot\Delta_q\,u\Vert_{L^p}^r\Big)^{\frac{1}{r}}\quad\text{si}\quad
r<\infty,\quad \text{et}\!\quad \Vert u\Vert_{ \dot B^s_{p,\infty}}
\overset{\rm d\acute{e}f}{=}\sup_{q\in\Z}\,2^{qs}\Vert\dot\Delta_q\,u
\Vert_{L^p}.
$$
 On d\'efinit
l'espace de Besov homog\`ene $\dot B^s_{p,r}:=\dot B^s_{p,r}(\R^N)$ par
\begin{itemize}

\item $\dot B^s_{p,r}=\{u\in\mathcal S'(\R^N)\; \; |\; \;
\|u\|_{\dot B^s_{p,r}}<+\infty\},$ lorsque $s<\frac{N}{p}$ ou
$s=\frac{N}{p}$ et $r=1$.
\item $\dot B^s_{p,r}=\{u\in\mathcal S'(\R^N)\; \; |\; \; \forall \;\;
|\alpha|=k\!+\!1,\;\; \partial^\alpha u\in \dot B^{s-k-1}_{p,r}\}$ si
$\frac{N}{p}+k\leq s<\frac{N}{p}+k+1,$  ou
$s=\frac{N}{p}+k+1$ et $r=1,$ pour un $k\in\N.$
\end{itemize}
\end{defi}
\begin{rem}
On prendra garde au fait que $\dot B^s_{p,r}(\R^N)$
est un espace de Banach si et seulement si 
$s<N/p$ ou $s\leq N/p$ et $r=1.$
\end{rem}
Rappelons \'egalement la d\'efinition
des espaces de Besov {\it non homog\`enes}~:
\begin{defi}
Soient $s\in\R,\,(p,r)\in[1,\infty]^2$ et $u\in {\cal {S}}'(\R^N).$
On note
$$
\begin{array}{lll}
\Vert u\Vert_{B^s_{p,r}}:=\Big(\Vert\dot S_0u\Vert_{L^p}^r+
\sum_{q\in\N}2^{rqs}
\Vert\dot\Delta_q\,u\Vert_{L^p}^r\Big)^{\frac{1}{r}}&\text{si}&
r<\infty,\\\Vert u\Vert_{B^s_{p,\infty}}
:=\max\Bigl(\Vert\dot S_0u\Vert_{L^p},\:\sup_{q\in\N}\,2^{qs}\Vert\dot\Delta_qu
\Vert_{L^p}\Bigr).&&\end{array}
$$
 L'espace de Besov non  homog\`ene
$B^s_{p,r}:=B^s_{p,r}(\R^N)$ est l'ensemble des distributions
tem\-p\'e\-r\'ees $u$ telles que $\Vert u\Vert_{B^s_{p,r}}$
soit fini.
\end{defi}
\begin{rem}
Les espaces de Besov
$B^s_{2,2}$ et $\dot B^s_{2,2}$ co\"{\i}ncident respectivement
avec les espaces de Sobolev $H^s$ et $\dot H^s.$
Si $s\in\R^+\setminus\N,$ les espaces de Besov $B^s_{\infty,\infty}$
et $\dot B^s_{\infty,\infty}$  co\"{\i}ncident respectivement
avec les espaces de  H\"older $C^s$ et  $\dot C^s.$
Par abus, nous noterons $C^s:=B^s_{\infty,\infty}$
et $\dot C^s:=\dot B^s_{\infty,\infty}$ pour tout $s\in\R.$
\end{rem}
Les in\'egalit\'es suivantes, dites {\it de Bernstein} 
et d\'emontr\'ees par exemple dans \cite{C} seront
d'un usage constant.
\begin{lem}  Soit $1\leq p_1\leq p_2\leq\infty$ et
$\psi\in\cC_c^\infty(\R^N).$  Alors on a
$$\hskip1cm\|\psi(2^{-q}D)u\|_{L^{p_2}}\leq C2^{qN\bigl(\frac 1{p_{1}}-\frac
1{p_{2}}\bigr)}\|\psi(2^{-q}D)u\|_{L^{p_1}}.$$
\end{lem}
 Comme cons\'equence de l'in\'egalit\'e de Bernstein  et
de la d\'efinition de $\dot B^s_{p,r},$ on a la proposition suivante~:
\begin{pro}\label{injection}
(i) Il existe une constante $c$ strictement positive telle que
\begin{equation}\label{equivalence}
c^{-1}\Vert u\Vert_{\dot B^s_{p,r}} \leq \Vert\nabla
u\Vert_{\dot B^{s-1}_{p,r}} \leq c \Vert u\Vert_{\dot B^s_{p,r}}.
\end{equation}
(ii) Pour $1\leq p_1\leq p_2\leq\infty$ et $1\leq r_1\leq r_2\leq\infty,$ on~a
$\dot B^s_{p_1,r_1}\hookrightarrow
\dot B^{s-N(\frac{1}{p_1}-\frac{1}{p_2})}_{p_2,r_2}.$\\
(iii) Si $p\in[1,\infty],$ alors
$\dot B^{\frac{N}{p}}_{p,1}\hookrightarrow
\dot B^{\frac{N}{p}}_{p,\infty}\cap L^{\infty}.$
Si $p$ est fini, l'espace $\dot B^{\frac Np}_{p,1}$ est
une alg\`ebre.\\
(iv) Interpolation r\'eelle: $(\dot B^{s_1}_{p,r},
\dot B^{s_2}_{p,r})_{\theta,\,r'}= \dot B^{\theta s_2+(1-\theta)s_1}_{p,r'}$
pour $0<\theta<1$ et $1\leq p,r,r'\!\leq~\!\!\infty.$
\end{pro}
Signalons aussi le r\'esultat d'inclusion suivant qui nous sera fort
utile~:
\begin{lem}\label{l:lorentz} Pour $1<p<q\leq\infty,$ on a
$$
L^{p,\infty}(\R^N)\hookrightarrow \dot B^{\frac Nq-\frac
Np}_{q,\infty}(\R^N).
$$
\end{lem}
\begin{dem}
Notons $h_\jmath=2^{\jmath N}h(2^\jmath\cdot)$ avec $h={\cal
F}^{-1}\varphi.$ On a $\dot\Delta_\jmath u=h_\jmath\star u,$ donc par les
in\'egalit\'es de convolution dans les espaces de
Lorentz  (voir par exemple \cite{L}),
$$
\|\dot\Delta_\jmath u\|_{L^q}\leq \|h_\jmath\|_{L^{r,1}}\|u\|_{L^{p,\infty}}
\quad\text{avec}\quad\frac1r=1+\frac1q-\frac1p\cdotp
$$
En faisant un changement de variable, on constate que
$$
\|h_\jmath\|_{L^{r,1}}=2^{\jmath N(1-\frac1r)}\|h\|_{L^{r,1}},
$$
d'o\`u le r\'esultat.
\end{dem}
Nous souhaitons maintenant donner quelques estimations a priori dans
les espaces de Besov pour l'\'equation de la chaleur. Ces
estimations d\'ecoulent du lemme suivant d\'emontr\'e dans~\cite{CH}~:
\begin{lem}\label{l:heatLP}
Il existe deux constantes $c$ et $C$ telles que pour tout
$\tau\geq0,$ $q\in\Z,$ et $p\in[1,\infty],$ on ait
$$
\|e^{\tau\Delta}\dot\Delta_qu\|_{L^p} \leq Ce^{-c\tau2^{2q}}\Vert
\dot\Delta_qu\Vert_{L^p}.
$$
\end{lem}
De ce lemme, on d\'eduit facilement (voir encore \cite{CH}) le
r\'esultat suivant~:
\begin{pro}\label{p:chaleur}
Soit $s\in\R,$ $1\leq p,r,\rho_1\leq\infty.$
Soit $u_0\in \dot B^s_{p,r}$
et $f\!\in\!\tilde L_T^{\rho_1}(\dot B^{s\!-\!2\!+\!\frac2{\rho_1}}_{p,r}).$
Alors l'\'equation  de la chaleur
$$
\d_tu-\nu\Delta u=f,\qquad u_{|t=0}=u_0
$$
admet une unique solution $u$ dans
$\tilde L_T^\infty(\dot B^s_{p,r})\cap\tilde L_T^{\rho_1}(\dot
B^{s+\frac2{\rho_1}}_{p,r})$
et  il existe une constante $C$ ne d\'ependant que de la dimension
$N$ telle que l'on ait l'estimation a priori suivante pour tout
$t\in[0,T]$ et $\rho\geq\rho_1$~:
\begin{equation}\label{eq:estchaleur}
\nu^{\frac1\rho}\Vert u\Vert_{\tilde L_t^\rho(\dot B^{s+\frac2\rho}_{p,r})}
\leq C\biggl(\Vert u_0\Vert_{\dot B^s_{p,r}}+
\nu^{\frac1{\rho_1}-1}\Vert f\Vert_{\tilde
L_t^{\rho_1}(\dot B^{s-2+\frac{2}{\rho_1}}_{p,r})}\biggr).
\end{equation}
Si $r<\infty,$ la solution $u$ appartient
 de plus \`a $\cC([0,T];\dot B^s_{p,r}).$
\end{pro}
Dans  la proposition ci-dessus, on a utilis\'e les
 espaces $\tilde L_T^\rho(\dot B^s_{p,r})$ qui se d\'efinissent
comme \`a la d\'efinition \ref{def} apr\`es avoir pos\'e
\begin{equation}\label{eq:ltilde}
\Vert u\Vert_{\tilde L_T^\rho(\dot B^s_{p,r})}:= \Bigl\|2^{qs}\Vert
\dot\Delta_qu\Vert_{L_T^\rho(L^p)}\Bigr\|_{\ell^r(\Z)}.
\end{equation}
Notons qu'en  vertu de l'in\'egalit\'e de Minkowski, on a
\begin{equation}\label{eq:Minkowski}
\Vert u\Vert_{L_T^\rho(\dot B^s_{p,r})}\leq
 \Vert u\Vert_{\tilde L_T^\rho(\dot B^s_{p,r})}
\quad\text{pour}\quad \rho\geq r,
\end{equation}
et l'in\'egalit\'e oppos\'ee si $\rho\leq r.$
\smallbreak
On posera $\tilde\cC_T(\dot B^s_{p,r}):=\tilde L_T^\infty(\dot B^s_{p,r})
\cap\cC([0,T];\dot B^s_{p,r})$ et
 $\tilde L^\rho_{loc}(\R_+;\dot B^s_{p,r})=\bigcap_{T>0}
\tilde L_T^\rho(\dot B^s_{p,r}).$

Sur le m\^eme mod\`ele, on peut d\'efinir
des espaces {\it non homog\`enes} $\tilde L_T^\rho(B^s_{p,r}).$
Dans le cas particulier o\`u $p=r=2,$ on utilisera plut\^ot
la notation $\tilde L_T^\rho(\dot H^s)$
ou  $\tilde L_T^\rho(H^s).$
\begin{rem}\label{r:stokes}
Gr\^ace \`a la proposition $\ref{p:chaleur}$, et en utilisant le fait
que le projecteur $\cP$ sur les champs \`a divergence nulle
est  un multiplicateur de Fourier homog\`ene de degr\'e
$0$ et est donc continu de $\dot B^s_{p,r}$ dans 
lui-m\^eme (voir e.g. \cite{L}), il est facile de r\'esoudre 
le syst\`eme de Stokes non stationnaire
\begin{equation}\label{eq:stokes}
\left\{\begin{array}{l}\d_tu-\nu\Delta u+\nabla\Pi=f,
\qquad{\rm div}\, u=0,\\
u_{|t=0}=u_0,\end{array} \right.
\end{equation}
avec donn\'ee initiale $u_0\in\dot B^s_{p,r}$ \`a divergence nulle, et
force $f\in\tilde L_T^1(\dot B^s_{p,r}).$

On obtient encore une unique solution $(u,\nabla\Pi)$ avec
$u\in\tilde L_T^\infty(\dot B^s_{p,r})\cap\tilde L_T^1(\dot
B^{s+2}_{p,r})$ et $\nabla\Pi\in\tilde
 L_T^1(\dot B^{s}_{p,r}),$
et $u$ v\'erifie pour tout $1\leq\rho\leq\infty$ l'estimation
$$\nu^{\frac1\rho}\Vert u\Vert_{\tilde L_T^\rho(\dot B^{s+\frac2\rho}_{p,r})}
\leq C\biggl(\Vert u_0\Vert_{\dot B^s_{p,r}} +\Vert \cP f\Vert_{\tilde
L_T^1(\dot B^{s}_{p,r})}\biggr).
$$
Si $r<\infty,$ la solution $u$ appartient
 de plus \`a $\cC([0,T];\dot B^s_{p,r}).$
\end{rem}
\begin{rem}\label{r:nonhomogene}
On peut  \'enoncer
une version non homog\`ene de la proposition
$\ref{p:chaleur}$ pour des donn\'ees
$u_0\in B^s_{p,r}$ et $f\in\tilde L_T^{\rho_1}(B^{s-2+\frac2{\rho_1}}_{p,r}).$
Les r\'esultats d'existence, d'unicit\'e et de continuit\'e
en temps demeurent, mais la constante $C$ de $\eqref{eq:estchaleur}$
se met \`a d\'ependre (au pire lin\'eairement) de $T.$
\end{rem}
 La preuve de certaines estimations a priori pour des termes
de type quadratique sera grandement
facilit\'ee par l'usage du  calcul
paradiff\'erentiel et notamment
 de la d\'ecomposition suivante introduite par J.-M. Bony
dans \cite{B}~:
$$
fg=\dot T_fg+\dot T_gf+\dot R(f,g),
$$
o\`u le terme de  paraproduit $\dot T$ est d\'efini par
$\displaystyle{\,\dot T_fg:=\sum_q\dot S_{q\!-\!1}f\,\dot\Delta_q g},$
et le terme
de reste $\dot R,$ par $\displaystyle{\,\dot R(f,g)
:=\sum_q \dot\Delta_qf\tilde\Delta_qg}$
avec $\tilde\Delta_q:=\dot\Delta_{q\!-\!1} +\dot\Delta_q
+\dot\Delta_{q\!+\!1}.$
\smallbreak
Nous \'enon\c cons ci-dessous
quelques r\'esultats de continuit\'e d'usage constant pour les
op\'erateurs de paraproduit et de reste
(consulter par exemple \cite{RS} pour une \'etude
exhaustive de ces op\'erateurs).
\begin{pro}\label{p:paraproduit}
Soit  $1\leq p,p_1,p_2,r,r_1,r_2\leq\infty$
v\'erifiant $\frac1p=\frac1{p_1}+\frac1{p_2}$ et
 $\frac1r=\frac1{r_1}+\frac1{r_2}.$
L'op\'erateur de paraproduit $\dot T$ est continu~:
\begin{itemize}
\item de $L^\infty\times\dot B^t_{p,r}$ dans
$\dot B^t_{p,r}$ pour tout $t\in\R,$
\item de $\dot B^{-s}_{p_1,r_1}\times
\dot B^t_{p_2,r_2}$ dans $\dot B^{t-s}_{p,r}$
pour tout $t\in\R$ et  $s>0.$
\end{itemize}
L'op\'erateur de reste $\dot R$ est continu~:
\begin{itemize}
\item de $\dot B^s_{p_1,r_1}\times\dot B^t_{p_2,r_2}$
dans $\dot B^{s+t}_{p,r}$ pour tout $(s,t)\in\R^2$ tel que $s+t>0,$
\item de $\dot B^s_{p_1,r_1}\times\dot B^{-s}_{p_2,r_2}$
dans $\dot B^{0}_{p,\infty}$ si $s\in\R$ et
 $\frac1{r_1}+\frac1{r_2}\geq1.$
\end{itemize}
\end{pro}
La proposition ci-dessus permet de montrer
la plupart des r\'esultats de continuit\'e
pour le produit de deux distributions
appartenant \`a des espaces de Besov.
Nous utiliserons fr\'equemment par exemple que
l'application $(u,v)\mapsto uv$ est continue de
\smallbreak
\begin{itemize}
\item
$\dot B^s_{N,1}\times\dot B^{t}_{N,1}$ dans
$\dot B^{s+t-1}_{N,1}$ si $-1<s,t\leq1$ et $s+t>0,$
\item $\dot H^s\times \dot H^t$ dans
$\dot B^{s+t-\frac N2}_{2,1}$
si $|s|<\frac N2,$  $|t|<\frac N2$ et $s+t>0.$
\end{itemize}


\section{Quelques estimations a priori}

\subsection{Espaces de Lorentz et syst\`eme
de Stokes non stationnaire}

Dans cette section, on d\'emontre quelques estimations
a priori dans les espaces de
Lorentz pour le syst\`eme de Stokes non stationnaire
\eqref{eq:stokes}.
Le lemme suivant nous permettra
de contr\^oler le terme de couplage avec l'\'equation
de transport dans le syst\`eme de Boussinesq.
 \begin{lem}\label{l:lorentz2}
On a les estimations a priori suivantes pour tout $t\geq0$~:
\begin{itemize}
\item Si $N=3\:$ alors
$\displaystyle{\; \nu\,\bigg\|\int_0^t e^{\nu(t-s)\Delta}\cP(\theta
e_3)(s)\,ds \bigg\|_{L^{3,\infty}}\leq
C\|\theta\|_{L^\infty_t(L^{1})}.}$
\item Si $N\geq4\:$ alors
$\displaystyle{\;\nu\,\bigg\|\int_0^t e^{\nu(t-s)\Delta} \cP(\theta
e_N)(s)\,ds\bigg\|_{L^{N,\infty}}\leq
C\|\theta\|_{L^\infty_t(L^{\frac N3,\infty})}}.$  \end{itemize}
\end{lem}
\begin{dem}
Par un changement de variable temporel, on peut se ramener au cas
$\nu=1.$ Notons par ailleurs
que $\cP\in \cL(L^p(\R^N);L^p(\R^N))$ pour tout $1<p<\infty.$
Par interpolation r\'eelle, on a donc en particulier
\begin{equation}\label{eq:P}
\cP\in\cL\bigl(L^{N,\infty}(\R^N),L^{N,\infty}(\R^N)\bigr).
\end{equation}
Par cons\'equent, on peut ``oublier" le projecteur $\cP$
dans les in\'egalit\'es \`a d\'emontrer et se limiter
\`a la preuve d'in\'egalit\'es similaires  pour l'op\'erateur de
la chaleur.
Pour ce faire, on \'ecrit la d\'ecomposition suivante avec
$A$ param\`etre positif~:
$$\displaylines{\int_0^t
e^{(t-s)\Delta}\theta(s)\,ds= \underbrace{\int_0^A
1_{[0,t]}(s)\,e^{(t-s)\Delta}\theta
(s)\,ds}_{I_A}+\underbrace{\int_A^\infty
1_{[0,t]}(s)\,e^{(t-s)\Delta}\theta(s)\,ds}_{I^A}.}$$ Par
calcul direct on a pour $\tau>0$ et $1\leq p\leq q\leq\infty,$
  \begin{eqnarray}\label{eq:G}
\|e^{\tau\Delta }\|_{{\cal L}(L^p,L^q)} =\|\mathcal
F^{-1}(e^{-\tau|\xi|^2})\|_{L^r}
=\frac{1}{\tau^{N/2}}\|G(\frac{x}{\sqrt\tau})\|_{L^r}
=\frac{1}{\tau^{N/2r'}}\|G\|_{L^r}\end{eqnarray} o\`u l'on a not\'e
$G=\mathcal F^{-1}(e^{-|\xi|^2}),$ $1/r=1+1/q-1/p$ et
$1/r'=1/p-1/q.$ \smallbreak

Supposons d'abord que  $N=3$. On a
$$\begin{array}{lll}\|I_A\|_{L^1}&\leq&
\Int_0^A\|e^{(t-s)\Delta}\theta(s)\|_{L^1}1_{[0,t]}(s)\,ds,\\[1.5ex]
 &\leq&
\Int_0^A\|e^{(t-s)\Delta}\|_{{\cal
L}(L^1,L^1)}\|\theta(s)\|_{L^1}1_{[0,t]}(s)\,ds,\\[1.5ex] &\leq&
\|G\|_{L^1}\Int_0^A\|\theta (s)\|_{L^1}1_{[0,t]}(s)\,ds\leq
A\|G\|_{L^1}\|\theta\|_{L^\infty_t(L^1)}.
\end{array}$$
D'autre part, pour $I^A$, on peut \'ecrire
$$
\begin{array}{lll}\|I^A\|_{L^\infty}&\leq&
\Int_A^\infty\|e^{(t-s)\Delta}\theta(s)\|_{L^\infty}1_{[0,t]}(s)\,ds,\\[2ex]
 &\leq&
\Int_A^\infty\|e^{(t-s)\Delta}\|_{{\cal
L}(L^1,L^\infty)}\|\theta(s)\|_{L^1}1_{[0,t]}(s)\,ds,\\[2ex]
&\leq&\|G\|_{L^\infty} \Int_A^\infty
\frac{1}{(t-s)^{3/2}}\|\theta (s)\|_{L^1}1_{[0,t]}(s)\,ds,\\[2ex]
&\leq&2A^{-1/2}\|G\|_{L^\infty}\|\theta\|_{L^\infty_t(L^1)}.
\end{array}
$$
En tenant compte de ces deux estimations,
de la remarque \ref{r:lorentz} et de \eqref{eq:P}, on obtient
le r\'esultat souhait\'e dans le cas $N=3.$
\smallbreak
Le cas $N\geq 4$  se traite de mani\`ere analogue. \`A partir de
 \eqref{eq:G},  on montre  par interpolation que
l'op\'erateur $e^{\tau\Delta}$ est continu de $L^{\frac N3,\infty}$ dans $L^{\frac
N3,\infty}$ avec norme born\'ee par $C\Vert G\Vert_{L^1}.$
Cela assure que pour tout $A\geq0,$ on~a
$$\begin{array}{lll}\|I_A\|_{L^{\frac N3,\infty}}&\leq&
\Int_0^A\|e^{(t-s)\Delta}\|_{{\cal L}(L^{\frac
N3,\infty},L^{\frac N3,\infty})}
\|\theta(s)\|_{L^{\frac N3,\infty}}1_{[0,t]}(s)\,ds,\\[1.5ex] &\leq&
CA\|G\|_{L^{1}}\|\theta\|_{L^\infty_t(L^{\frac N3,\infty})}.
\end{array}$$
Pour majorer $I^A$, on utilise le fait  que
l'op\'erateur  $e^{\tau\Delta}$ est continu de $L^{\frac
N3,\infty}$ dans $L^{\infty}$ avec norme
$C\tau^{-\frac32}\Vert G\Vert_{L^{\frac N{N-3},1}}$
(combiner  \eqref{eq:G} avec un argument d'interpolation).
 On peut donc  \'ecrire
$$
\begin{array}{lll}\|I^A\|_{L^\infty} &\leq&
\Int_A^\infty\|e^{(t-s)\Delta}\|_{{\cal L}(L^{\frac
N3,\infty},L^\infty)}
\|\theta(s)\|_{L^{\frac N3,\infty}}1_{[0,t]}(s)\,ds,\\[2ex]
&\leq&C\|G\|_{L^{\frac N{N-3},1}} \Int_A^\infty
\frac{1}{(t-s)^{3/2}}\|\theta (s)\|_{L^{\frac N3,\infty}}1_{[0,t]}(s)\,ds,\\[2ex]
&\leq&CA^{-1/2}\|G\|_{L^{\frac N{N-3},1}}
\|\theta\|_{L^\infty_t(L^{\frac N3,\infty})}.
\end{array}
$$
En combinant les majorations pour $I_A$ et $I^A,$  le fait que
$$
(L^{\frac N3,\infty},L^\infty)_{\frac23,\infty}=L^{N,\infty},
$$
et \eqref{eq:P}, on peut
maintenant conclure \`a l'in\'egalit\'e souhait\'ee dans le cas
$N\geq4.$
\end{dem}
Le lemme suivant
(qui semble \^etre d\^u \`a Y. Meyer, voir
 \cite{Meyer})
nous permettra de contr\^oler le terme
de convection dans l'\'equation sur~$u.$
\begin{lem}\label{l:lorentz1}
Pour tout $t\geq0$ et $\nu\geq0,$
on a les estimations a priori
$$\nu\,\bigg\|\int_0^te^{\nu(t-s)\Delta}\cP
\nabla g(s)\,ds\bigg\|_{L^{N,\infty}}\leq
\left\{\begin{array}{lll}
C\|g\|_{L_t^\infty(L^{1})}&\text{ si }&N=2,\\
C\|g\|_{L_t^\infty(L^{\frac N2,\infty})}&\text{ si }&N\geq3.\end{array}\right.
$$
 \end{lem}
\begin{dem}
Pour la commodit\'e du lecteur, nous donnons ci-dessous les grandes lignes de
la d\'emonstration. On proc\`ede comme dans le lemme pr\'ec\'edent~:
apr\`es s'\^etre ramen\'e au cas $\nu=1,$
on d\'ecoupe l'int\'egrale \`a majorer en deux parties
$J_A$ et $J^A.$
En vertu de \eqref{eq:P}, on peut ``oublier" le projecteur
$\cP.$
De plus, il est clair que pour tout $\tau>0$ et $1\leq p\leq q\leq\infty,$ on a
\begin{equation}\label{eq:H}
\Vert e^{\tau\Delta}\nabla\Vert_{\cL(L^p,L^q)}
=\frac1{\tau^{\frac12+\frac N{2r'}}}\Vert H\Vert_{L^r}
\quad\!\text{avec}\!\quad H:={\cal F}^{-1}\bigl(e^{-|\eta|^2}\eta\bigr)
\ \text{ et }\ \frac1{r'}=\frac1p-\frac1q\cdotp
\end{equation}
En reprenant les calculs du lemme pr\'ec\'edent,
on obtient donc si $N\geq3,$
$$
\Vert J_A\Vert_{L^{\frac N2,\infty}}\leq CA\Vert H\Vert_{L^1}\Vert
g\Vert_{L^{\frac N2,\infty}}\quad \text{et}\quad
\Vert{J^A}\Vert_{L^\infty}\leq
CA^{-1}\Vert H\Vert_{L^{\frac N{N-2},1}}\Vert
g\Vert_{L^{\frac N2,\infty}}.
$$
On conclut alors en remarquant que
$(L^{\frac N2,\infty},L^\infty)_{\frac12,\infty}=L^{N,\infty}.$

Le cas $N=2$ est formellement identique. Il suffit
 de remplacer $L^{\frac N2,\infty}$ par $L^1.$
\end{dem}


\subsection{Estimations a priori pour le syst\`eme
de Navier-Stokes}

Dans cette partie, on d\'emontre des estimations
a priori pour  le syst\`eme de Navier-Stokes incompressible~:
\begin{equation}
\label{Navier-Stokes}
\begin{cases} \partial_t u+u\cdot\nabla
u-\nu\Delta u+\nabla\Pi=f,\\
{\rm div}\, u=0\\
u(0,x)=u_0(x)
\end{cases}
\end{equation}
sous diverses conditions de petitesse.
\smallbreak
\'Enon\c cons tout d'abord un r\'esultat global
de propagation  de la r\'egularit\'e
$L^\infty_t(\dot B^{-1+\frac Np}_{p,1})\cap
L^1_t(\dot B^{1+\frac Np}_{p,1})$ sous l'hypoth\`ese que $u$ est
petite dans $\dot C^{-1}.$
\begin{lem}\label{l:stokes}
Soit $p\in[1,\infty[$ et  $u$ une solution   de
$\eqref{Navier-Stokes}$
appartenant \`a $L^\infty_T(\dot B^{-1+\frac Np}_{p,1})\cap
L^1_T(\dot B^{1+\frac Np}_{p,1}).$
Il existe deux  constantes $c$ et $C$  telles que si $u_0$
appartient \`a $\dot B^{-1+\frac Np}_{p,1}$ et v\'erifie
 ${\rm div}\, u_0=0,$  si
 $f$ est dans $L^1([0,T];\dot B^{-1+\frac Np}_{p,1})$
et de plus
$$\sup\limits_{t\in [0,T]}\|u(t)\|_{\dot C^{-1}}\leq c\nu$$
alors, on a pour
tout $t\in [0,T]$ l'estimation
$$\|u(t)\|_{\dot B^{\frac Np-1}_{p,1}}+\nu\int_0^t
\|u(\tau)\|_{\dot B^{\frac Np+1}_{p,1}}\,d\tau \leq C\bigg(\|u_0\|_{\dot B^{\frac
Np-1}_{p,1}}+\int_0^t\|f(\tau)\|_{\dot B^{\frac Np-1}_{p,1}}\,d\tau\bigg).$$
\end{lem}

\begin{dem} On applique l'op\'erateur
$\cP\dot\Delta_q $ \`a l'\'equation de Navier-Stokes. Sachant que $\div
u=0,$ on obtient
$$
\d_t\dot\Delta_qu-\nu\Delta\dot\Delta_qu=\cP\dot\Delta_qf
-\cP\dot\Delta_q(u\cdot\nabla u).
$$
Admettons un instant que
\begin{equation}\label{eq:convection}
\|\cP\dot\Delta_q(u\cdot\nabla u)\|_{L^p}\lesssim a_q 2^{-q(-1+\frac
Np)}\|u\|_{\dot C^{-1}}\|u\|_{\dot B^{1+\frac Np}_{p,1}}\quad\text{avec}\quad
\sum_{q\in\Z}a_q=1.
\end{equation}
Alors  en utilisant le lemme \ref{l:heatLP}, on a
$$\displaylines{
2^{q(-1+\frac Np)}\|\dot\Delta_q u(t)\|_{L^p}\leq 2^{q(-1+\frac
Np)}e^{-c\nu t2^{2q}}\|\dot\Delta_q u_0\|_{L^p}\hfill\cr\hfill +\int_0^t
e^{-c\nu(t-\tau)2^{2q}}\Bigl(2^{q(-1+\frac Np)}\|\ddq f\|_{L^p}+
a_q\|u\|_{\dot C^{-1}}\|u\|_{\dot B^{-1+\frac Np}_{p,1}}\Bigr)\,d\tau.}$$
Apr\`es
 sommation sur $q,$ on trouve donc
$$\displaylines{\|u\|_{L^\infty_t(\dot B^{\frac Np-1}_{p,1})}
+\nu\|u\|_{L^1_t(\dot B^{\frac Np+1}_{p,1})}
\leq C\Big(\|u_0\|_{\dot B^{-1+\frac
Np}_{p,1}}\hfill\cr\hfill+\|f\|_{L^1_t(\dot B^{-1+\frac Np}_{p,1})}+
\Vert u\Vert_{L^\infty_t(\dot C^{-1})}\Vert u\Vert_{L^1_t(\dot B^{1+\frac
Np}_{p,1})}\Big).}$$ Si l'on suppose que
$$C\|u\|_{L^\infty_t(\dot C^{-1})}\leq \frac\nu 2$$
alors on obtient  l'estimation souhait\'ee.
\end{dem}
L'in\'egalit\'e \eqref{eq:convection} que nous avons utilis\'ee plus
haut repose sur le lemme suivant~: \begin{lem}\label{l:convection}
Soient $s>1$, $r\in [1,+\infty]$ et $u$ un champ \`a divergence
nulle. Alors il existe une constante $C$ et une suite
$a_q\in\ell^r(\Z)$ telles que
$$\|\dot\Delta_q(u\cdot\nabla v)\|_{L^p}\leq C2^{-q(s-2)}a_q
\bigl(\|u\|_{\dot C^{-1}}\|v\|_{\dot B^s_{p,r}}
+\|v\|_{\dot C^{-1}}\|u\|_{\dot B^s_{p,r}}\bigr)\quad\text{et}\quad
\Vert (a_q)\Vert_{\ell^r}\leq1.$$
\end{lem}
\begin{dem}
En vertu de \eqref{presortho} et de $\div u=0,$ on a
$$
\dot\Delta_q(u\cdot\nabla v)= \sum\limits_{q'\geq
q-4}\dot\Delta_q\div(\dot S_{q'+2}u\otimes\dot\Delta_{q'}v)
+\sum\limits_{|q-q'|\leq 4}\dot\Delta_q(\dot\Delta_{q'} u\cdot\nabla
\dot S_{q'-1}v).$$
 Sachant que  $s>1,$ on a
$$
\begin{array}{lll}
\sum\limits_{q'\geq q-4}\|\dot\Delta_q\div(\dot S_{q'+2}u\otimes\dot\Delta_{q'}v)\|_{L^p}&\leq&
C2^q\sum\limits_{q'\geq q-4}\|\dot S_{q'+2}u\|_{L^\infty}\|\dot\Delta_{q'}v\|_{L^p},\\
&\leq& C2^q\sum\limits_{q'\geq q-4}2^{-q'(s-1)}
\|u\|_{\dot C^{-1}}\,2^{q's}\|\dot\Delta_{q'}v\|_{L^p},\\
&\leq& C2^{-q(s-2)}a_q\|u\|_{\dot C^{-1}}\|v\|_{\dot B^s_{p,r}}
\end{array}
$$
avec $a_q\in\ell^r(\Z)$ v\'erifiant $\Vert (a_q)\Vert_{\ell^r}\leq1,$
\smallbreak
D'autre part, on dispose de la majoration suivante~:
$$
\begin{array}{lll}
\sum\limits_{|q'-q|\leq 4}\|\dot\Delta_q(\dot\Delta_{q'}u\cdot\nabla
\dot S_{q'-1}v)\|_{L^p}&\leq& C\sum\limits_{|q'-q|\leq 4}
\|\dot\Delta_{q'}u\|_{L^p}\|\dot S_{q'-1}\nabla v\|_{L^\infty},\\
&\leq& C2^{-q(s-2)}a_q\|v\|_{\dot C^{-1}}\|u\|_{\dot B^s_{p,r}},
\end{array}
$$
d'o\`u le lemme.
\end{dem}
Il n'est pas clair que le lemme \ref{l:stokes}
demeure valable dans le cas limite $p=+\infty.$
Mais on dispose toutefois du substitut suivant~:
\begin{lem}\label{l:stokeslimit}
Soit $u$ une solution   de $\eqref{Navier-Stokes}$
 avec  donn\'ee initiale  $u_0\in\dot B^{-1}_{\infty,1}$
\`a divergence nulle, et terme de force  $f\in L^1([0,T], \dot
B^{-1}_{\infty,1}).$
Il existe deux  constantes $c$ et $C$  telles que si
$$\sup\limits_{t\in [0,T]}\|u(t)\|_{L^{N,\infty}}\leq c\nu$$
alors on a  pour tout
$t\in [0,T]$ l'estimation
$$\|u(t)\|_{\dot B^{-1}_{\infty,1}}+\nu\|u\|_{L^1_t(\dot
B^{1}_{\infty,1})} \leq C\bigg(\|u_0\|_{\dot B^{-1}_{\infty,1}}
+\|f\|_{L^1_t(\dot B^{-1}_{\infty,1})}\bigg).$$
 \end{lem}
\begin{dem}
La d\'emarche g\'en\'erale est la m\^eme que dans la preuve du lemme
\ref{l:stokes}, mais on ne peut plus utiliser le lemme \ref{l:convection} qui
est  faux dans le cas limite $p=\infty,$ $r=1$ et $s=1.$
Le bon substitut est l'in\'egalit\'e suivante~:
\begin{equation}\label{eq:bonsubstitut}
\|\dot\Delta_q\bigl(u\cdot\nabla u\bigr)\|_{L^\infty}\leq Ca_q2^q
\|u\|_{L^{N,\infty}}\|u\|_{\dot B^1_{\infty,1}}\quad\text{avec}\quad
\sum_{q\in\Z}a_q\leq1.
 \end{equation}
La justification de \eqref{eq:bonsubstitut}
repose sur la d\'ecomposition de
Bony. En tenant compte de $\div u=0,$ on a
pour tout $j\in\{1,\cdots,N\},$
$$\bigl[\div(u\otimes u)\bigr]^j=\d_i(\dot T_{u^i}u^j+\dot T_{u^j}u^i+\dot R(u^i,u^j)).$$
 Les r\'esultats de continuit\'e pour le paraproduit
\'enonc\'es dans la proposition \ref{p:paraproduit} assurent que
$$\|\div\dot\Delta_q (\dot T_u u)\|_{L^\infty}\leq C
2^q a_q\|u\|_{\dot C^{-1}}\|u\|_{\dot B^1_{\infty,1}} \quad\text{avec}\quad
\sum_{q\in\Z}a_q\leq1$$
et l'on a vu que $L^{N,\infty}\hookrightarrow\dot C^{-1}$ (cf lemme
\ref{l:lorentz}).
\smallbreak
La continuit\'e de l'injection de $L^{N,\infty}$
dans $\dot C^{-1}$ permet aussi d'\'ecrire que
$$\|\dot\Delta_q \div \dot R(u,u)\|_{L^\infty}\leq C2^{q}\|\dot\Delta_q \div \dot R(u,u)\|_{L^{N,\infty}}\quad\text{pour tout}\quad
q\in\Z.$$
Par cons\'equent, on~a
$$\begin{array}{lll}
\|\dot\Delta_q \div \dot R(u,u)\|_{L^\infty}
&\!\!\leq\!\!&
C 2^{2q}\,\Bigl\|\sum\limits_{q'\geq q-3}\dot\Delta_q(\dot\Delta_{q'}u
\tilde\Delta_{q'}u)\Bigr\|_{L^{N,\infty}},\\
&\!\!\leq\!\!&
C2^{2q}\sum\limits_{q'\geq q-3}\|\dot\Delta_{q'}u\|_{L^{N,\infty}}
\|\tilde\Delta_{q'} u\|_{L^\infty},\\
&\!\!\leq\!\!& C  2^{q}\sum\limits_{q'\geq
q-3}2^{q-q'}a_{q'}\|u\|_{L^{N,\infty}}\|u\|_{\dot B^1_{\infty,1}}
\end{array}$$
avec $\sum_qa_q\leq1,$
ce qui ach\`eve la d\'emonstration de \eqref{eq:bonsubstitut}.
\end{dem}


\subsection{Estimations pour l'\'equation de transport}

Dans cette section nous \'etablissons diverses
estimations a priori pour l'\'equation de transport
\begin{equation}
\label{transport}
\begin{cases}
\partial_t\theta+\div(\theta u)=f\\
\theta(0,x)=\theta_0(x)
\end{cases}
\end{equation}
avec $u$  champ de vecteur lipschitizien {\it \`a  divergence
nulle}.
Pour commencer, rappelons un r\'esultat classique
relatif aux espaces
$L^p$ ou $L^{p,\infty}$:
\begin{pro}\label{p:transportLp}
Pour tout $p\in[1,\infty],$ on a les estimations a priori
suivantes~:
$$
\forall t\in\R^+,\;\|\theta (t)\|_{L^p}\leq
\|\theta_0\|_{L^p}+\int_0^t\|f(s)\|_{L^p}\,ds.$$ Plus
g\'en\'eralement, pour $1\leq p<\infty,$  on a
$$\forall t\in\R^+,\;\|\theta (t)\|_{L^{p,\infty}}\leq
\|\theta_0\|_{L^{p,\infty}}+\int_0^t\|f(s)\|_{L^{p,\infty}}\,ds.$$
Dans tous les cas, l'in\'egalit\'e peut \^etre remplac\'ee par une
\'egalit\'e si $f\equiv0.$
\end{pro}
\begin{dem}
En introduisant  le flot $\psi$  du champ $u$ d\'efini par
\begin{equation}
\begin{cases}
\partial_t\psi(t,x)=u(t,\psi(t,x)),\\
\psi(0,x)=x,
\end{cases}
\end{equation}
on obtient \begin{equation}\label{eq:flot}
\theta(t,x)=\theta_0(\psi(t,x)^{-1})+\int_0^tf(s,\psi(s,\psi(t,x)^{-1}))\,ds.
\end{equation}
Sachant que $\div u=0$, la fonction $\psi_t:=\psi(t,\cdot)$ est un
diff\'eomorphisme qui pr\'eserve le volume.
En cons\'equence,  $$\|\theta (t)\|_{L^p}\leq
\|\theta_0\|_{L^p}+\int_0^t\|f(s)\|_{L^p}$$
avec \'egalit\'e si $f\equiv 0.$
\smallbreak
 La preuve de
l'in\'egalit\'e dans les espaces $L^{p,\infty}$ se traite en
comparant les ensembles de niveau de $\theta$ et de $\theta_0.$ En
effet~: dans le cas $f\equiv0,$ on a pour tout $\lambda\geq0,$
\begin{equation}\label{eq:niveau}
\bigl\{y\in\R^N\,/\, |\theta(t,y)|>\lambda\bigr\}
=\psi_t\Bigl(\bigl\{x\in\R^N\,/\,
|\theta_0(x)|>\lambda\bigr\}\Bigr),
\end{equation}
et les deux ensembles ont m\^eme mesure puisque $\psi$ pr\'eserve la
mesure.

Le cas g\'en\'eral $f\not\equiv0$ s'en d\'eduit facilement \`a
partir de la formule \eqref{eq:flot}.
\end{dem}
Dans la proposition ci-dessous, nous
\'enon\c cons un r\'esultat de  propagation
de r\'egularit\'e de type Besov pour les \'equations de
transport.
\begin{pro}\label{p:transport1}
Soient $1\leq p\leq p_1\leq+\infty,$ $r\in [1,\infty]$ et
 $s\in\R$ tels que
$$
-1-\min\Bigl(\frac N{p_1},\frac N{p'}\Bigr) <s<1+\frac
N{p_1}\cdotp$$ Supposons que le champ de vecteurs $u$ soit \`a
divergence nulle et \`a coefficients dans
$L^1([0,T];\dot B^{1+\frac{N}{p_1}}_{p_1,1}(\R^N)),$
 que le  terme de force $f$ appartienne \`a
$L^1([0,T];\dot B^{s}_{p,r}(\R^N))$ et que la donn\'ee initiale
$\theta_0$ soit dans $\dot B^{s}_{p,r}(\R^N)$. Alors l'\'equation
$(\ref{transport})$ admet une unique solution $\theta\in
L^\infty([0,T];\dot B^{s}_{p,r})$ ($\theta\in\cC([0,T];\dot B^s_{p,r})$ si
$r<\infty$)
 qui v\'erifie l'estimation
$$\|\theta\|_{\tilde L_t^\infty(\dot B^{s}_{p,r})}
\leq\Bigl(\|\theta_0\|_{\dot B^{s}_{p,r}}+
\|f\|_{\tilde L_t^1(\dot B^{s}_{p,r})}\Bigr)
\exp\biggl(C\int_0^t\|\nabla u\|_{\dot B^{\frac
N{p_1}}_{p_1,1}}d\tau\biggr).
$$
L'in\'egalit\'e ci-dessus reste valable dans le cas limite
$r=\infty$ et
 $s=-1-\min\Bigl(\frac N{p_1},\frac N{p'}\Bigr).$
\end{pro}
\begin{dem}
Un r\'esultat analogue
a \'et\'e \'etabli dans \cite{DANCHI}
(mais dans le cadre des espaces de Besov {\it non homog\`enes}
et sans le cas limite).
Nous nous limitons donc
 \`a une  preuve succincte des estimations a priori.

Appliquons $\dot\Delta_q $ \`a l'\'equation de transport \eqref{transport}.
 En utilisant  la d\'ecomposition de Bony, la propri\'et\'e
de divergence nulle et \eqref{presortho}, on obtient
$$
\d_t\dot\Delta_q\theta+\dot S_{q-1}u\cdot\nabla\dot\Delta_q\theta
=\dot\Delta_qf-F_q,
$$
$$
\displaylines{\text{avec}\quad F_q=\dot\Delta_q (\dot T_{\partial_i\theta}u^i)
+\d_i\dot\Delta_q\dot R(u^i,\theta)\hfill\cr\hfill+\sum_{|q-q'|\leq
1}(\dot S_{q'\!-\!1}\!-\!\dot S_{q\!-\!1})u\cdot\nabla\dot\Delta_q\dot\Delta_{q'}\theta
+\sum\limits_{|q-q'|\leq 4}[\dot\Delta_q,
\dot S_{q'-1}u]\cdot\nabla\dot\Delta_{q'}\theta.}
$$
D'apr\`es la proposition \ref{p:transportLp}, on a
\begin{equation}\label{eq:theta}
\Vert\dot\Delta_q \theta(t)\Vert_{L^{p}}\leq \Vert\dot\Delta_q
\theta_0\Vert_{L^{p}}+ \int_0^t \Vert\dot\Delta_q
f(\tau)\Vert_{L^{p}}\,d\tau +\int_0^t \Vert
F_q(\tau)\Vert_{L^{p}}\,d\tau.
\end{equation}
Il s'agit maintenant de majorer $\Vert F_q\Vert_{L^{p}}$
de fa\c con ad\'equate.
Le premier terme se traite en faisant appel \`a
la proposition \ref{p:paraproduit}.
Comme $s<1+\frac N{p_1},$ on obtient
$$\Vert\dot\Delta_q\dot T_{\d_i\theta}u^i\Vert_{L^{p}}
\leq C 2^{-qs}a_q \Vert\nabla
u\Vert_{\dot B^{\frac{N}{p_1}}_{p_1,1}}\|\theta\|_{\dot B^{s}_{p,r}}
\quad\text{avec}\quad \Vert a_q\Vert_{\ell^r}\leq 1.$$
De m\^eme,
sous la condition $s+1+\min\Bigl(\frac N{p_1},\frac N{p'}\Bigr)>0,$
la proposition \ref{p:paraproduit} donne
\begin{equation}\label{eq:vit2}
\Vert\d_i\dot\Delta_q\dot R(u^i,\theta)\Vert_{L^{p}} \leq C a_q
2^{-qs}\Vert\nabla u\Vert_{\dot B^{\frac{N}{p_1}}_{{p_1},1}}
\Vert\theta\Vert_{\dot B^{s}_{p,r}}\quad\text{avec}\quad
\Vert a_q\Vert_{\ell^r}\leq 1,
\end{equation}
et  dans le cas limite  $s+1+\Bigl(\frac N{p_1},\frac N{p'}\Bigr)=0,$
\begin{equation}\label{eq:vit3}
\sup_{q\in\Z}\Vert\d_i\dot\Delta_q\dot R(u^i,\theta) \Vert_{L^{p}} \leq
C2^{-qs}\Vert\nabla
u\Vert_{\dot B^{\frac{N}{p_1}}_{{p_1},1}}
\Vert\theta\Vert_{\dot B^{s}_{p,\infty}}.
\end{equation}
 Le  troisi\`eme terme de $F_q$ est une fonction localis\'ee spectralement
dans une couronne dyadique $2^q\mathcal C.$ Ainsi, sans aucune
restriction sur $s,$ on~a
$$
\sum_{\vert
q-q'\vert\leq1}\Vert(\dot S_{q\!-\!1}-\dot S_{q'\!-\!1})u\!\cdot\!\nabla\dot\Delta_q\dot\Delta_{q'}\theta\Vert_{L^{p}}
\leq Ca_q2^{-qs} \Vert\nabla u
\Vert_{\dot B^{\frac N{p_1}}_{p_1,1}}\Vert\theta\Vert_{\dot B^s_{p,r}}.
$$
Enfin, par une
in\'egalit\'e classique sur les commutateurs (voir par exemple
\cite{C}),
on~a
$$
\Vert[\dot\Delta_q, \dot S_{q'-1}u]\cdot\nabla\dot\Delta_{q'}\theta\Vert_{L^{p}}
\lesssim
 2^{-q}\Vert\nabla \dot S_{q'-1}u\Vert_{L^\infty}
\Vert\dot\Delta_{q'}\nabla \theta\Vert_{L^{p}}.
$$
En cons\'equence,
$$
\Vert[\dot\Delta_q, \dot S_{q'-1}u]\cdot\nabla\dot\Delta_{q'}\theta\Vert_{L^{p}}
\lesssim 2^{q'-q}\Vert\nabla
u\Vert_{L^\infty}\Vert\dot\Delta_{q'}\theta\Vert_{L^{p}},
$$
d'o\`u
$$
\sum\limits_{|q-q'|\leq 4}\Vert[\dot\Delta_q,
\dot S_{q'-1}u]\cdot\nabla\dot\Delta_{q'}\theta\Vert_{L^{p}} \leq C
2^{-qs}a_q\Vert\nabla u\Vert_{L^\infty}
\Vert\theta\Vert_{\dot B^{s}_{p,r}}\quad\text{avec}\quad
\Vert(a_q)\Vert_{\ell^r}\leq1.
$$
En sommant toutes ces estimations et en utilisant le fait que
$\dot B^{\frac N{p_1}}_{p_1,1}\hookrightarrow L^\infty,$ on trouve
(en supposant  $r=\infty$ dans le cas limite)~:
$$\Vert
F_q\Vert_{L^{p}} \lesssim 2^{-qs}a_q \Vert\nabla
u\Vert_{\dot B^{\frac{N}{p_1}}_{p_1,1}}\Vert \theta\Vert_{\dot B^{s}_{p,r}}
\quad\text{avec}\quad \Vert a_q\Vert_{\ell^r}\leq 1.
$$
 Il ne reste plus qu'\`a injecter cette in\'egalit\'e dans
\eqref{eq:theta},
\`a multiplier par $2^{qs}$ puis
\`a prendre la norme $\ell^r(\Z)$ en $q$ de chaque membre. On obtient
$$
\Vert\theta\Vert_{\tilde L_t^\infty(\dot B^{s}_{p,r})} \leq
\Vert\theta_0\Vert_{\dot B^{s}_{p,r}} +\Vert
f\Vert_{\tilde L_t^1(\dot B^{s}_{p,r})}+C\int_0^t\Vert\nabla
u(\tau)\Vert_{\dot B^{\frac{N}{p_1}}_{p_1,1}}
\Vert\theta(\tau)\Vert_{\dot B^{s}_{p,r}}d\tau
$$
et le lemme de Gronwall
permet alors d'obtenir l'estimation souhait\'ee.
\end{dem}
La pr\'esence du terme exponentiel dans l'estimation
a priori de la proposition \ref{p:transport1} semble \^oter tout espoir
de d\'emontrer des r\'esultats globaux pour le
syst\`eme de Boussinesq.

Dans \cite{Vishik} cependant, M. Vishik a remarqu\'e que
cette croissance exponentielle n'a pas lieu (dans le cas  $f\equiv0$)
si l'on travaille dans  des espaces de Besov \`a {\it indice
de r\'egularit\'e nul}.
 Le r\'esultat de M. Vishik
a \'et\'e red\'emontr\'e r\'ecemment par
T. Hmidi et S. Keraani
par une m\'ethode plus robuste qui
permet entre autres d'ajouter un terme de diffusion
  dans l'\'equation de transport (voir \cite{HK}).
Nous utiliserons la version suivante du
 r\'esultat de Vishik~:
\begin{pro}\label{p:transport2}
Soit $\theta\in\tilde L^\infty_T(\dot B^0_{p,r})$ une solution de
$\eqref{transport}$ avec  $f\in\tilde L^1_T(\dot B^0_{p,r})$ et
 $\nabla u\in L^1([0,T], \dot B^{0}_{\infty,1})$  \`a divergence
nulle.
Alors on a la majoration suivante pour tout $t\in[0,T]$:
$$\|\theta\|_{\tilde L_t^\infty(\dot B^0_{p,r})}\leq C
\Bigl(\|\theta_0\|_{\dot B^0_{p,r}} +\|f\|_{\tilde L^1_t(\dot B^0_{p,r})}\Bigr)
\Bigl(1+\int_0^t\|\nabla u(s)\|_{\dot B^{0}_{\infty,1}}\,ds\Bigr).$$
\end{pro}
\begin{dem} Nous adoptons l'approche utilis\'ee dans
\cite{HK}.
Pour simplifier la pr\'esen\-tation, on se limite au cas $r=1$
(qui est le seul utilis\'e dans l'article).

Tout d'abord, remarquons que
par unicit\'e de la solution de l'\'equation de transport, on a
$\theta=\sum\limits_{q\in\Z} \theta_q$ o\`u $\theta_q$
est la solution de l'\'equation
\begin{equation}
\begin{cases}
\partial_t\theta_q+u\cdot\nabla\theta_q=\dot\Delta_q f\\
\theta_q(0,x)=\dot\Delta_q\theta_0(x).
\end{cases}
\end{equation}
 On a clairement
\begin{equation}\label{eq:besov0}
\|\theta\|_{\dot B^0_{p,1}}\leq\sum\limits_{j,q\in\Z}\|\dot\Delta_j
\theta_q\|_{L^p}=\sum\limits_{|j-q|\leq Q}\|\dot\Delta_j\theta_q
\|_{L^p}+\sum\limits_{|j-q|\geq Q}\|\dot\Delta_j\theta_q \|_{L^p},
\end{equation}
o\`u $Q$ est un param\`etre entier que l'on ajustera plus loin.
\smallbreak
 En faisant appel \`a la
proposition \ref{p:transportLp}, on constate que
$$\begin{array}{lll}
\sum\limits_{|j-q|\leq Q}\|\dot\Delta_j\theta_q(t)
\|_{L^p}&\lesssim&\sum\limits_{|j-q|\leq Q}\|\theta_q(t)\|_{L^p},\\
&\lesssim& Q\sum\limits_{q}\bigl(\|\dot\Delta_q\theta_0
\|_{L^p}+\|\dot\Delta_q f\|_{L^1_t(L^p)}\bigr),\\&\lesssim&
Q(\|\theta_0\|_{\dot B^0_{p,1}}+\|f\|_{L^1_t(\dot B^0_{p,1})}).
\end{array}
$$
Pour traiter la deuxi\`eme somme de \eqref{eq:besov0},
on utilise la proposition
\ref{p:transport1}
 avec donn\'ee initiale $\ddq\theta_0$
dans $\dot B^{\pm \epsilon}_{p,1}$
et terme de force dans $L^1(0,T;\dot B^{\pm\epsilon}_{p,1})$
(avec $\epsilon$ choisi dans $(0,1)$).
 On peut donc \'ecrire
$$\|\theta_q(t)\|_{\dot B^{\pm \epsilon}_{p,1}}
\leq (\|\dot\Delta_q \theta_0\|_{\dot B^{\pm \epsilon}_{p,1}}+\|\dot\Delta_q
f\|_{L^1_t(\dot B^{\pm \epsilon}_{p,1})})\exp\bigg(C\int_0^t\|\nabla
u\|_{\dot B^{0}_{\infty,1}}\,d\tau\bigg),$$
  ce qui implique
$$\|\dot\Delta_j\theta_q(t)\|_{L^p}
\leq 2^{-\epsilon|q-j|}a_j(\|\dot\Delta_q \theta_0\|_{L^p}+\|\dot\Delta_q
f\|_{L^1_t(L^p)})\exp\bigg(C\int_0^t\|\nabla
 u\|_{\dot B^{0}_{\infty,1}}\,d\tau\bigg),$$
avec $\Vert(a_j)\Vert_{\ell^1(\Z)}\leq1.$ Par sommation sur les indices
$(j,q)$ tels que $|j-q|\geq Q,$ on trouve
$$\sum\limits_{|j-q|\geq Q}\|\dot\Delta_j\theta_q(t)\|_{L^p}
\leq 2^{-Q\epsilon}C(\| \theta_0\|_{\dot B^{0}_{p,1}}+\|
f\|_{L^1_t(\dot B^{0}_{p,1})})\exp\bigg(C\int_0^t\|\nabla
u\|_{\dot B^{0}_{\infty,1}}\,d\tau\bigg).$$
Il ne reste plus qu'\`a choisir $Q$ tel que
$Q\epsilon\log2\approx\Bigl(1+C\Int_0^t\|\nabla u\|_{\dot
  B^{0}_{\infty,1}}\,d\tau\Bigr).$
 \end{dem}


\subsection{Estimations \`a pertes}

La d\'emonstration de l'unicit\'e
 pour le syst\`eme de Boussinesq en dimension deux
repose sur  des estimations  a priori
avec perte de r\'egularit\'e pour les \'equations de
transport-diffusion du type
\begin{equation}\label{eq:td}
\begin{cases}
\partial_t \rho+\div(\rho u)-\nu\Delta\rho=f\\
\rho|_{t=0}=\rho_0
\end{cases}
\end{equation}
 avec un champ de vecteurs $u$
\`a divergence nulle qui  n'est pas lipschitzien en espace.

Dans le cas sans diffusion,  H. Bahouri et J.-Y. Chemin
ont remarqu\'e (voir \cite{BCh}) que si le champ
de vitesses \'etait  localement int\'egrable en temps
\`a valeurs log-lipschitz, i.e
$$|u(t,x)-u(t,y)|\leq\gamma(t)|x-y|(1-\log|x-y|)\quad\text{pour tout
}\hskip0.3cm|x-y|<1\quad\text{avec}\quad\gamma\in L^1_{loc}(0,T)$$
alors  la r\'egularit\'e de la donn\'ee
initiale \'etait  susceptible de se d\'egrader au cours du temps. Par exemple,
si  $\rho_0$ appartient \`a l'espace de Sobolev $H^s$ alors la solution \`a
l'instant $t$  est  seulement  dans
 $H^{s-V(t)}$ avec $V(t)=C\int_0^t\gamma(\tau)\,d\tau.$

Nous souhaitons d\'emontrer des estimations
analogues valables pour tout $\nu\geq0$,
sous des hypoth\`eses l\'eg\`erement plus  g\'en\'erales que celles
de \cite{BCh}  et
avec en sus un gain de r\'egularit\'e de type parabolique
si $\nu>0.$
\smallbreak
Nous avons opt\'e
pour un cadre fonctionnel {\it non homog\`ene}.
De ce fait, nous ferons appel \`a
des blocs dyadiques  {\it non homog\`enes} d\'efinis comme suit~:
$$
\dq:=\ddq \quad\text{si}\quad q\geq0,\qquad
\Delta_{-1}:=\dot S_0\quad\text{et}\quad \dq:=0\quad\text{si}\quad q\leq-2,
$$
et \`a la troncature basses fr\'equences
$$
S_q:=\sum_{p\leq q-1}\Delta_p.
$$
\begin{pro}\label{loglip} Soit $s\in]-1-\frac N2,1+\frac N2[$
et $\rho$ une solution de l'\'equation
de transport-diffusion $\eqref{eq:td}$.
Il existe $N_0\in\N$ ne d\'ependant que du support de
la fonction $\varphi$ intervenant dans la d\'ecomposition
de Littlewood-Paley, une constante universelle $C_0$
et  deux constantes $c$ et $C$ ne d\'ependant que de $s$ et de $N$
  tels que si
\begin{equation}\label{eq:condloglip}
\Vert u\Vert_{\tilde L_T^1(B^{1+\frac N2}_{2,\infty})}\leq c
\end{equation}
alors on a l'estimation a priori suivante pour
tout $t\in[0,T]$:
$$\displaylines{
\sup_{\substack{q\geq-1\\
\tau\in[0,t]}}2^{qs-\e_q(\tau)}\Vert\dq\rho(\tau)\Vert_{L^2}
+\nu\sup_{q\geq0}\int_0^t2^{q(s+2)-\e_q(\tau)}
\Vert\dq\rho(\tau)\Vert_{L^2}\,d\tau
\hfill\cr\hfill\leq C_0\biggl(\Vert\rho_0\Vert_{B^s_{2,\infty}}
+\sup_{q\geq-1}\int_0^t2^{qs-\e_q(\tau)}
\Vert\dq f(\tau)\Vert_{L^2}\,d\tau
\biggr)}
$$
o\`u $\e_q(t):=C\Sum_{q'=-1}^q2^{q'(1+\frac
N2)}\Int_0^t\Vert\check\Delta_{q'}u\Vert_{L^2}\,d\tau$
avec $\check\Delta_{q'}=\sum_{|\alpha|\leq N_0}\Delta_{q'+\alpha}$
v\'erifie
 $$\e_q(T)-\e_{q'}(T)\leq\frac12\bigl(1+\frac
N2+s\bigr)(q-q')\quad\text{pour}\quad\!\! q\geq q'.$$
\end{pro}
\begin{dem}
 La d\'emonstration de cette proposition se trouve implicitement dans
\cite{DANCHI} (voir th. 3.12) mais l'effet
r\'egularisant du Laplacien n'y est pas explicit\'e.
\'Etant donn\'e que l'exploitation de cet effet r\'egularisant
est la clef de la preuve de
l'unicit\'e pour le syst\`eme de Boussinesq
en dimension $2,$ nous donnons ici la d\'emonstration int\'egrale
de l'estimation a priori.
\smallbreak
En appliquant l'op\'erateur $\Delta_q$ \`a \eqref{eq:td}
et en proc\'edant comme dans la preuve de la proposition \ref{p:transport1},
on obtient
$$
\d_t\dq\rho+S_{q-1}u\cdot\nabla\dq\rho-\nu\Delta\dq\rho=\dq f
+F_q^1+F_q^2+F_q^3+F_q^4
$$
avec
$$
\begin{array}{lll}
F_q^1:=\Sum_{|q'-q|\leq4}[S_{q'-1}u,\dq]\cdot\nabla\Delta_{q'}\rho,&&
F_q^2:=\Sum_{|q'-q|\leq1}\bigl(S_{q-1}-S_{q'-1}\bigr)u\cdot
\nabla\dq\Delta_{q'}\rho,\\
F_q^3:=-\dq\Bigl(\Sum_{|q'\!-\!q|\leq4}\!\!S_{q'-1}\d_i\rho\,
\Delta_{q'}u^i\Bigr),&&
F_q^4:=-\!\Sum_{q'\geq q-3}\!\d_i\dq\biggl(
\Delta_{q'}\rho\Bigl(\Sum_{|\alpha|\leq1}\Delta_{q'\!+\!\alpha}\Bigr)
u^i\biggr).
\end{array}
$$
Une m\'ethode d'\'energie standard combin\'ee
avec l'in\'egalit\'e de Bernstein donne
\begin{equation}\label{eq:unicite1}
\frac12\frac d{dt}\Vert\dq\rho\Vert_{L^2}^2
+\nu_q2^{2q}\Vert\dq\rho\Vert_{L^2}^2
\leq \Bigl(\Vert\dq f\Vert_{L^2}
+\sum_{i=1}^4\Vert F_q^i\Vert_{L^2}\Bigr)\Vert\dq\rho\Vert_{L^2}
\end{equation}
avec $\nu_q:=\kappa\nu$ si $q\geq0$ et $\nu_{-1}=0,$  la constante
$\kappa$ d\'ependant uniquement du support de $\varphi.$
\smallbreak
Par des manipulations classiques, on~a
$$
\begin{array}{ll}
\!\!\Vert F_q^1\Vert_{L^2}\leq\!\Sum_{|q'\!-\!q|\leq4}\!
\Vert{\nabla S_{q'\!-\!1}u}\Vert_{L^\infty}
\Vert\Delta_{q'}\rho\Vert_{L^2},&
\!\Vert F_q^2\Vert_{L^2}\leq\Sum_{|q'-q|\leq1}
2^q\Vert{\check\Delta_q u}\Vert_{L^\infty}\Vert\dq\rho\Vert_{L^2},\\
\!\!\Vert F_q^3\Vert_{L^2}\leq\!\Sum_{q'\leq q+2}\!2^{q'}
\Vert{\Delta_{q'}\rho}\Vert_{L^\infty}
\Vert{\check\Delta_q u}\Vert_{L^2},&
\!\Vert F_q^4\Vert_{L^2}\leq\!\!\Sum_{q'\geq q-3}\!2^{q(1\!+\!\frac N2)}
\Vert\Delta_{q'}\rho\Vert_{L^2}\Vert\check\Delta_{q'}u\Vert_{L^2},
\end{array}
$$
en notant $\check\Delta_q=\sum_{|\alpha|\leq N_0}\Delta_{q'}$ avec
$N_0$ suffisamment grand.
\smallbreak
En utilisant le fait que
$$\displaylines{
\Vert{\nabla S_{q'-1}u}\Vert_{L^\infty}\lesssim
\sum_{q''<q'-N_0}2^{q''(1+\frac N2)}\Vert\check\Delta_{q''}u\Vert_{L^2},\cr
\Vert{\check\Delta_q u}\Vert_{L^\infty}\lesssim
2^{q\frac N2}\Vert{\check\Delta_qu}\Vert_{L^2},\quad
\Vert{\Delta_{q'}\rho}\Vert_{L^\infty}\lesssim2^{q'\frac N2}
\Vert{\Delta_{q'}\rho}\Vert_{L^2},}
$$
et en injectant les majorations obtenues ci-dessus dans
\eqref{eq:unicite1}, on obtient
$$\displaylines{\frac12\frac d{dt}\Vert\dq\rho\Vert_{L^2}^2
+\nu_q2^{2q}\Vert\dq\rho\Vert_{L^2}^2
\leq\Bigl( \Vert\dq f\Vert_{L^2}
\hfill\cr\hfill
+C\sum_{q'\geq q-4}2^{(q-q')(1+\frac N2)}
\e'_{q'}\Vert\Delta_{q'}\rho\Vert_{L^2}
+C\sum_{q'\leq q}2^{(q'-q)(1+\frac N2)}
\e'_{q}\Vert\Delta_{q'}\rho\Vert_{L^2}\Bigr)\Vert\dq\rho\Vert_{L^2}}
$$
avec $\e_p(t):=\Int_0^t\sum_{p'\leq p}2^{p'(1+\frac N2)}
\Vert\check\Delta_pu\Vert_{L^2}\,d\tau.$
\smallbreak
Pour $\lambda>0,$ on pose
$$\rho_q^\lambda(t):=2^{qs}e^{-\lambda\e_q(t)}\Vert\dq\rho(t)\Vert_{L^2}\quad
\text{et}\quad
f_q^\lambda(t):=2^{qs}e^{-\lambda\e_q(t)}\Vert\dq f(t)\Vert_{L^2}
$$
de telle sorte que
$$\displaylines{
\frac12\frac d{dt}(\rho_q^\lambda)^2+\lambda\e'_q(\rho_q^\lambda)^2
+\nu_q2^{2q}(\rho_q^\lambda)^2\leq\rho_q^\lambda
\biggl(f_q^\lambda\hfill\cr\hfill+C2^{qs}e^{-\lambda\e_q}
\Bigl(\sum_{q'\geq q-4}2^{(q-q')(1+\frac N2)}
\e'_{q'}\Vert\Delta_{q'}\rho\Vert_{L^2}
+\sum_{q'\leq q}2^{(q'-q)(1+\frac N2)}
\e'_{q}\Vert\Delta_{q'}\rho\Vert_{L^2}\Bigr)\biggr).}
$$
De l'in\'egalit\'e pr\'ec\'edente, on d\'eduit apr\`es int\'egration
en temps~:
$$
\displaylines{\rho_q^\lambda(t)
+\nu_q2^{2q}\int_0^t\rho_q^\lambda(\tau)\,d\tau
+\lambda\int_0^t\e'_q(\tau)\rho_q^\lambda(\tau)\,d\tau
\leq\rho_q^\lambda(0)+\int_0^tf_q^\lambda(\tau)\,d\tau\hfill\cr\hfill
+C\sum_{q'\geq q-4}2^{(q-q')(1+\frac N2+s)}
\int_0^t\e'_{q'}(\tau)e^{\lambda(\e_{q'}(\tau)-\e_q(\tau))}
\rho_{q'}^\lambda(\tau)\,d\tau\hfill\cr\hfill
+C\sum_{q'\leq q}2^{(q'-q)(1+\frac N2-s)}\int_0^t
\e'_{q}(\tau)e^{\lambda(\e_{q'}(\tau)-\e_q(\tau))}
\rho_{q'}^\lambda(\tau)\,d\tau.}
$$
En \'ecrivant $\e'_q(\tau)=\e'_{q'}(\tau)
+(\e'_q-\e'_{q'})(\tau)$ dans le dernier terme
et en remarquant que
 $(\e_n)_{n\geq-1}$ est une
suite croissante de  fonctions croissantes positives,
on obtient pour tout $q\geq-1,$
$$
\displaylines{\rho_q^\lambda(t)+\nu_q2^{2q}\int_0^t\rho_q^\lambda(\tau)\,d\tau
+\lambda\int_0^t\e'_q(\tau)\rho_q^\lambda(\tau)\,d\tau
\leq\rho_q^\lambda(0)+\int_0^tf_q^\lambda(\tau)\,d\tau
\hfill\cr\hfill
+C\sum_{q'\geq q}2^{(q-q')(1+\frac N2+s)}e^{\lambda(\e_{q'}(t)-\e_q(t))}
\int_0^t\e'_{q'}(\tau)\rho_{q'}^\lambda(\tau)\,d\tau\hfill\cr\hfill
+C\sum_{q'\leq q}2^{(q'-q)(1+\frac N2-s)}\int_0^t
\e'_{q'}(\tau)\rho_{q'}^\lambda(\tau)\,d\tau+\frac C\lambda
\sum_{q'\leq q}2^{(q'-q)(1+\frac N2-s)}\sup_{\tau\in[0,t]}
\rho_{q'}^\lambda(\tau).}
$$
Supposons d\'esormais  que
\begin{equation}\label{eq:unicite2}
\sup_{q\geq-1}2^{q\bigl(1+\frac N2\bigr)}\Vert\check\Delta_q
u\Vert_{L^1_t(L^2)}\leq\e\log2
\quad\text{pour un }\ \e\ \text{ tel que }\
\lambda\e\leq\frac12\bigl(\frac N2+1+s\bigr)
\end{equation}
de telle sorte que pour tout $q'\geq q\geq-1,$ on ait
$$
2^{(q-q')(1+\frac N2+s)}e^{\lambda(\e_{q'}(t)-\e_q(t))}
\leq 2^{\frac{q-q'}2(1+\frac N2+s)}.
$$
Alors en prenant la borne sup\'erieure en $q$
de l'in\'egalit\'e pr\'ec\'edente, on obtient
$$\displaylines{
\sup_{\substack{q\geq-1\\
\tau\in[0,t]}}\biggl(\rho_q^\lambda(\tau)
+\nu_q2^{2q}\int_0^\tau\rho_q^\lambda(\tau')\,d\tau'
+\lambda\int_0^\tau\e'_q(\tau')\rho_q^\lambda(\tau')\,d\tau'\biggr)
\leq\sup_{q\geq-1}\rho_q^\lambda(0)\hfill\cr\hfill
+\sup_{q\geq-1}\int_0^tf_q^\lambda(\tau)\,d\tau
+C\sup_{q\geq-1}\int_0^t\e'_{q}(\tau)\rho_{q}^\lambda(\tau)\,d\tau
+\frac C\lambda\sup_{\substack{q\geq-1\\
\tau\in[0,t]}}\rho_q^\lambda(\tau).}
$$
En cons\'equence, on~a,
$$\displaylines{
\sup_{\substack{q\geq-1\\
\tau\in[0,t]}}\rho_q^\lambda(\tau)
+\nu\sup_{q\geq0}2^{2q}\int_0^t\rho_q^\lambda(\tau)\,d\tau
+\lambda\sup_{q\geq-1}\int_0^t\e'_q(\tau)\rho_q^\lambda(\tau)\,d\tau
\leq3\Vert\rho_0\Vert_{B^s_{2,\infty}}\hfill\cr\hfill
+3\sup_{q\geq-1}\int_0^tf_q^\lambda(\tau)\,d\tau
+3C\sup_{q\geq-1}\int_0^t\e'_{q'}(\tau)\rho_{q'}^\lambda(\tau)\,d\tau
+\frac{3C}\lambda\sup_{\substack{q\geq-1\\
\tau\in[0,t]}}\rho_q^\lambda(\tau).}
$$
Il ne reste plus qu'\`a choisir $\lambda=6C$ pour conclure \`a
l'in\'egalit\'e souhait\'ee~:
$$
\sup_{\substack{q\geq-1\\
\tau\in[0,t]}}\rho_q^\lambda(\tau)
+\nu\sup_{q\geq0}2^{2q}\int_0^t\rho_q^\lambda(\tau)\,d\tau
\leq C_0\biggl(\Vert\rho_0\Vert_{B^s_{2,\infty}}
+\sup_{q\geq-1}\int_0^tf_q^\lambda(\tau)\,d\tau\biggr).
$$
\end{dem}
\begin{rem}\label{r:loglip}
L'estimation de la proposition \ref{loglip}
s'\'etend sans aucune difficult\'e au probl\`eme
de Stokes non stationnaire avec convection~:
$$
\left\{\begin{array}{l}
\d_t\rho+{\rm div}\,(\rho u)-\nu\Delta\rho+\nabla\Pi=f,\\
{\rm div}\, u=0.\end{array}\right.
$$
La propri\'et\'e de divergence nulle assure en effet
que le terme $\nabla\Pi$ est ``transparent'' et
que les estimations sur $\Vert\dq\rho\Vert_{L^2}^2$ restent
les m\^emes que pour l'\'equation $\eqref{eq:td}$.
\end{rem}


\section{Solutions fortes}

Cette partie est consacr\'ee \`a la d\'emonstration
de r\'esultats d'existence locale (ou globale 
\`a donn\'ees petites) et d'unicit\'e pour
le syst\`eme de Boussinesq.

\subsection{Un r\'esultat local}

Dans un premier temps, nous souhaitons 
\'etablir le  th\'eor\`eme \ref{th:quatre}.
Pour cela, nous allons adopter le sch\'ema classique suivant~:
\begin{enumerate}
\item Construction d'une
suite de  solutions globales approch\'ees,
\item Estimations uniformes en grande norme sur un (petit) intervalle
de temps fixe,
\item Convergence en petite norme sur cet intervalle de temps,
\item Unicit\'e.
\end{enumerate}
Le crit\`ere de prolongement sera d\'emontr\'e \`a la fin de la section.
\begin{enumerate}
\item \underline{Construction d'une
suite de  solutions globales approch\'ees.}

Notons $(u_L,\nabla\Pi_L)$ la solution libre de l'\'equation de
Stokes
$$
\d_tu_L-\nu\Delta u_L+\nabla\Pi_L=0,\qquad u_L(0)=u_0.
$$
On a bien s\^ur (cf remarque \ref{r:stokes})~:
$$
u_L\in\cC_b(\R_+;\dot B^0_{N,1})\cap L^1(\R_+;\dot B^2_{N,1})\quad\text{et}
\quad\nabla\Pi_L\in L^1(\R_+;\dot B^0_{N,1}).
$$
On d\'efinit ensuite  $(\theta^0,\ov u^0)\equiv(\theta_0,0)$ puis on
r\'esout par r\'ecurrence les \'equations lin\'eaires suivantes~:
\begin{equation}\label{schema}
\begin{cases}
\partial_t\theta^{n+1}+u^n\cdot\nabla\theta^{n+1}=0\quad\text{avec}\quad
u^n:=u_L+\ov u^n,\\
\d_t\ov u^{n+1}-\nu\Delta\ov u^{n+1}+\nabla\ov\Pi^{n+1}
=-\div(u^n\otimes u^n)+\theta^{n+1}e_N,\\
\theta^{n+1}(0)=\theta_0,\qquad u^{n+1}(0)=0.
\end{cases}
\end{equation}
D'apr\`es la proposition \ref{p:transport1}, la remarque
\ref{r:stokes}  et le fait que $\dot B^1_{N,1}$ soit une
alg\`ebre, on obtient ainsi une suite
$(\theta^n,\ov u^n,\nabla\ov\Pi^n)_{n\in\N}$ de solutions globales
approch\'ees v\'erifiant
$$
\theta^n\in\cC(\R_+;\dot B^0_{N,1}),\quad \ov
u^n\in\cC(\R_+;\dot B^0_{N,1})\cap
L^1_{loc}(\R_+;\dot B^2_{N,1})\quad\text{et} \quad\nabla\ov \Pi^n\in
L^1_{loc}(\R_+;\dot B^0_{N,1}).
$$
\item \underline{Estimations a priori uniformes  en grande norme
sur un intervalle de temps fixe.}

D'apr\`es la proposition \ref{p:transport2}, on a pour tout
$t\geq0,$
\begin{equation}\label{eq:exist1}
\Vert\theta^{n+1}(t)\Vert_{\dot B^0_{N,1}} \leq
C\Vert\theta_0\Vert_{\dot B^0_{N,1}} \Bigl(1+\Vert
u^n\Vert_{L_t^1(\dot B^2_{N,1})}\Bigr).
\end{equation}
Par ailleurs, d'apr\`es la remarque \ref{r:stokes} et les lois de produit
dans les espaces de Besov \'enonc\'ees apr\`es la proposition
\ref{p:paraproduit}, on~a
$$\displaylines{
\ov U^{n+1}(t):=\Vert\ov u^{n+1}\Vert_{L_t^1(\dot B^0_{N,1})}
+\nu\Vert\ov u^{n+1}\Vert_{L_t^1(\dot B^2_{N,1})}
+\Vert\nabla\ov\Pi^{n+1}\Vert_{L_t^1(\dot B^0_{N,1})}
\hfill\cr\hfill\lesssim\Vert u^{n}\Vert_{L_t^2(\dot B^1_{N,1})}^2
+\Vert\theta ^{n+1}\Vert_{L_t^1(\dot B^0_{N,1})}.}
$$
D'o\`u, en ins\'erant l'in\'egalit\'e \eqref{eq:exist1} dans
l'in\'egalit\'e ci-dessus,
$$\displaylines{
\ov U^{n+1}(t) \leq C\nu^{-1}\ov U^n(t)\bigl(\ov
U^n(t)
+t\Vert\theta_0\Vert_{\dot B^0_{N,1}}\bigr)\hfill\cr\hfill
+C\Bigl(\Vert u_L\Vert_{L_t^2(\dot B^1_{N,1})}^2
+t\Vert\theta_0\Vert_{\dot B^0_{N,1}}
\bigl(1+\Vert u_L\Vert_{L_t^1(\dot B^2_{N,1})}\bigr)\Bigr)}
$$
avec $C\geq1$ ne d\'ependant que de la dimension $N.$ \smallbreak
On en d\'eduit que  $\ov U^n(t)\leq\nu/(2C)$
tant que $t$ v\'erifie
\begin{equation}\label{eq:temps}
4C^2\nu^{-1}\Bigl(\Vert u_L\Vert_{L_t^2(\dot B^1_{N,1})}^2
+t\Vert\theta_0\Vert_{\dot B^0_{N,1}}
(2+\Vert u_L\Vert_{L_t^1(\dot B^2_{N,1})})\Bigr)\leq1,
\end{equation}
 En revenant \`a \eqref{eq:exist1} et en notant
$T$ la borne sup\'erieure des temps $t$ v\'erifiant
\eqref{eq:temps},  on peut finalement conclure que
$(\theta^n,u^n,\nabla\Pi^n)_{n\in\N}$ est  uniform\'ement born\'ee
dans  l'espace $E_T.$

\item \underline{Convergence en petite norme.}

Soit $\ov\theta^n:=\theta^n-\theta_0.$
Nous allons montrer que, quitte \`a diminuer $T$ si n\'ecessaire,
 la suite $(\ov\theta^n,\ov u^n,\nabla\ov\Pi^n)_{n\in\N}$
 est de Cauchy dans l'espace
$$F_T:=
\cC([0,T];\dot B^{-1}_{N,1})\times\Bigl(\cC([0,T];\dot B^{-1}_{N,1}) \times
L^1([0,T];\dot B^1_{N,1})\Bigr)^N\times
\Bigl(L^1([0,T];\dot B^{-1}_{N,1})\Bigr)^N.
$$
Comme les espaces de Besov utilis\'es sont homog\`enes,
l'appartenance \`a $E_T$
 n'entra\^{\i}ne a priori pas que
 $(\ov\theta^n,\ov u^n,\nabla\ov\Pi^n)$ est dans $F_T.$
En s'aidant de
 \eqref{schema} et de la proposition \ref{p:paraproduit},
 il n'est cependant pas difficile de v\'erifier que
$\d_t\ov\theta_n$ et  $\d_t\ov u_n$ sont dans
$L^2(0,T;\dot B^{-1}_{N,1}).$ Sachant que $\ov\theta_n(0)$
et $\ov u_n(0)$ sont nuls, on peut alors facilement justifier
que $(\ov\theta^n,\ov u^n,\nabla\ov\Pi^n)$ est dans $F_T.$
\smallbreak
Notons $\dt^{n}:=\theta^{n\!+\!1}\!-\!\theta^n,$
 $\du^n:=u^{n+1}-u^n$ et $\dPi^n:=\Pi^{n+1}-\Pi^n.$
Comme $\dt^n$ v\'erifie
$$
\d_t\dt^{n}+u^{n-1}\cdot\nabla\dt^n=-\div(\du^{n-1}\theta^n),
$$
la proposition \ref{p:transport1} assure  que
$$\Vert\dt^n(t)\Vert_{\dot B^{-1}_{N,1}}
\leq Ce^{\Vert u^{n}\Vert_{L_t^1(\dot B^2_{N,1})}}
\int_0^t\Vert\theta^n\,\du^{n\!-\!1}\Vert_{\dot B^0_{N,1}}\,d\tau.
$$
Donc, gr\^ace aux bornes uniformes obtenues dans l'\'etape
pr\'ec\'edente, on a
\begin{equation}\label{eq:exist3}
\Vert\dt^n(t)\Vert_{\dot B^{-1}_{N,1}} \leq
C_T\int_0^t\Vert\du^{n\!-\!1}\Vert_{\dot B^1_{N,1}}\,d\tau,
\end{equation}
avec $C_T$ ind\'ependante de $n$ et de $t\in[0,T].$ \smallbreak
Ensuite, comme
$$
\d_t\du^n-\nu\Delta\du^n+\nabla\dPi^n=
-\div(u^{n\!-\!1}\otimes\du^{n\!-\!1})
-\div(\du^{n\!-\!1}\otimes u^n)+\dt^n e_N,
$$
on a, d'apr\`es la remarque \ref{r:stokes},
$$\displaylines{
\dU^n(t):= \Vert\du^n\|_{L_t^\infty(\dot B^{-1}_{N,1})}
+\nu\Vert\du^n\|_{L_t^1(\dot B^{1}_{N,1})}
+\Vert\nabla\dPi^n\|_{L_t^1(\dot B^{-1}_{N,1})} \hfill\cr\hfill\leq
C\int_0^t\Bigl( \Vert
u^{n-1}\otimes\du^{n\!-\!1}\Vert_{\dot B^{0}_{N,1}}+
\Vert\du^{n-1}\otimes u^n\Vert_{\dot B^{0}_{N,1}}
+\Vert\dt^n\Vert_{\dot B^{-1}_{N,1}}\Bigr)\,d\tau.}
$$
Le dernier terme  peut \^etre born\'e gr\^ace \`a \eqref{eq:exist3}.
En utilisant  les lois de produit dans les
espaces de Besov \'enonc\'ees apr\`es la proposition
\ref{p:paraproduit}, et les in\'egalit\'es d'interpolation
de la proposition \ref{injection}, on obtient
$$\displaylines{
\Vert u^{n-1}\otimes\du^{n\!-\!1}\Vert_{L_t^1(\dot B^{0}_{N,1})}+
\Vert\du^{n-1}\otimes u^n\Vert_{L_t^1(\dot B^{0}_{N,1})}
\hfill\cr\hfill\leq C\bigl(\Vert u^{n-1}\Vert_{L_t^2(\dot B^{1}_{N,1})} +
\Vert u^{n}\Vert_{L_t^2(\dot B^{1}_{N,1})}\bigr)
\Vert\du^{n-1}\Vert_{L_t^2(\dot B^{0}_{N,1})}.}
$$
En combinant avec \eqref{eq:exist3}, on a donc finalement
$$ \dU^n(t)\leq \dU^{n-1}(t)
\nu^{-1/2}\Bigl(C_T\,t
+C\bigl(\Vert u^{n-1}\Vert_{L_t^2(\dot B^{1}_{N,1})}
+ \Vert u^{n}\Vert_{L_t^2(\dot B^{1}_{N,1})}\bigr)\Bigr).
$$
 D'apr\`es \eqref{eq:temps},
 il existe une constante $c$ (que l'on peut choisir
aussi petite que l'on veut quitte \`a diminuer $T$)
telle~que
$$
\Vert u^{n-1}\Vert_{L_T^2(\dot B^{1}_{N,1})}
+ \Vert u^{n}\Vert_{L_T^2(\dot B^{1}_{N,1})}\leq c\sqrt\nu.
$$
On en conclut que pour $T$ suffisamment petit, on~a
$$
\forall t\in[0,T],\;\forall
n\in\N,\;\dU^n(t)\leq\frac12\dU^{n-1}(t).
$$
La suite consid\'er\'ee est donc de
Cauchy dans $F_T.$

 \item \underline{Fin de la preuve de l'existence locale.}

On a donc d'une part $(\ov\theta^n,\ov u^n,\nabla\ov\Pi^n)$ qui converge dans
l'espace $F_T$
vers une fonction $(\ov\theta,\ov u,\nabla\ov\Pi),$ et d'autre part les
estimations uniformes en grande norme de la premi\`ere \'etape. Cela
permet d'affirmer que $(\theta^n,u^n,\nabla\Pi^n)_{n\in\N}$ converge vers une
solution $$\bigl(\theta:=\ov\theta+\theta_0,u:=\ov
u+u_L,\nabla\Pi:=\nabla(\Pi_L+\ov\Pi)\bigr)$$ du syst\`eme de
Boussinesq appartenant~\`a
\begin{equation}\label{eq:convergence}
L^\infty_T(\dot B^0_{N,1}) \times\Bigl(L^\infty_T(\dot B^0_{N,1})\cap
L^1_T(\dot B^2_{N,1})\Bigr)^N \times\Bigl(L^1_T(\dot B^0_{N,1})\Bigr)^N.
\end{equation}
La continuit\'e en temps s'obtient \`a l'aide
de la proposition \ref{p:transport1} et de la remarque
\ref{r:stokes} en remarquant que
\eqref{eq:convergence} entra\^ine que $\d_t\theta+\div(\theta u)$ et
$\d_tu-\nu\Delta u+\nabla\Pi$ sont dans $L^1_T(\dot B^0_{N,1}).$

\item\underline{Unicit\'e}

L'unicit\'e
s'obtient par les m\^emes arguments que la convergence en petite
norme. Les d\'etails sont laiss\'es au lecteur.

\item\underline{Crit\`ere de prolongement.}

Sous les hypoth\`eses du th\'eor\`emes,
on constate que $\d_t\theta+u\cdot\nabla\theta$ et $\d_tu-\nu\Delta
u+\nabla\Pi$ sont dans $L^1(0,T^*;\dot B^0_{N,1}).$
D'apr\`es la proposition \ref{p:transport1} et la remarque
\ref{r:stokes}, les fonctions $\theta$ et $u$ sont donc en fait dans
$\tilde L_T^\infty(\dot B^0_{N,1}).$ En partant de n'importe que temps
$t\in[0,T^*),$ on peut alors, gr\^ace \`a \eqref{eq:temps} obtenir
l'existence d'une solution avec donn\'ees $(\theta(t),u(t))$ sur un
intervalle de temps {\it ind\'ependant} de $t.$ Combin\'e avec la
propri\'et\'e d'unicit\'e,
 cela permet de prolonger $(u,\theta)$
au-del\`a de $T^*.$
\end{enumerate}
Cela ach\`eve la d\'emonstration du th\'eor\`eme \ref{th:quatre}.
\hfill\rule{2.1mm}{2.1mm}


\subsection{R\'esultats globaux \`a donn\'ees petites}

Cette section est consacr\'ee \`a la d\'emonstration
du th\'eor\`eme \ref{th:deux}.
Pour all\'eger la pr\'esen\-tation, on convient d\'esormais
que $L^{\frac N3,\infty}$ d\'esigne $L^1$ si $N=3.$

\subsubsection{Unicit\'e}

La partie unicit\'e du th\'eor\`eme  \ref{th:deux}
d\'ecoule de la proposition suivante~:
\begin{pro}\label{p:uniqueness}
Soit $(\theta_1,u_1,\nabla\Pi_1)$ et $(\theta_2,u_2,\nabla\Pi_2)$
 deux solutions du
syst\`eme de Boussinesq telles que $(\theta_1,u_1)$
et $(\theta_2,u_2)$
appartiennent
\`a l'espace
 $$\Bigl({\cal
C}([0,T];\dot B^0_{N,1})\cap L^\infty(0,T;L^{\frac
N3,\infty})\Bigl)\times
\Bigl({\cal C}([0,T];\dot B^{-1}_{\infty,1})\cap L^1(0,T;\dot B^1_{\infty,1})
\cap L^\infty(0,T;L^{N,\infty})\Bigr)^N. $$
Si les deux solutions  co\"{\i}ncident \`a l'instant initial alors elles
co\"{\i}ncident sur $[0,T].$
\end{pro}
\begin{p}
Remarquons tout d'abord que du fait des hypoth\`eses,
 on~a
$$\d_tu_i-\mu\Delta u_i\in L_T^1(\dot
B^{-1}_{\infty,1})\quad\text{pour}\quad i=1,2.$$
 La proposition  \ref{p:chaleur} assure donc que
$u_i\in\tilde L_T^\infty(\dot B^{-1}_{\infty,1}).$
Ce fait sera utilis\'e \`a plusieurs reprises dans les calculs qui suivent.
\smallbreak
Soit $(\dt,\du,\dPi)=(\theta_1-\theta_2,u_1-u_2,\Pi_1-\Pi_2).$
Remarquons que
\begin{equation}
\label{unicite}
\begin{cases}
\partial_t\dt+\div(u_2\,\dt)=-\div(\du\, \theta_1),
\\ \partial_t\du+\div(u_2\otimes\du)+\div(\du\otimes u_1)-\nu\Delta \du
+\nabla\dPi=\dt\,e_N.
\end{cases}
\end{equation}
Admettons provisoirement que $(\dt,\du,\nabla\dPi)$ appartienne  \`a l'espace
$$
G_T:=L^\infty\big([0,T];\,\dot B^{-1}_{N,\infty}\big) \times
\bigg(\cC\big([0,T];\,\dot B^{-1}_{N,\infty}\big)\cap \widetilde
L^1_T\big(\dot B^{1}_{N,\infty}\big)\bigg)^N \times\biggl(\widetilde
L^1_T\big(\dot B^{-1}_{N,\infty}\big)\biggr)^N.
$$
On munit $G_T$ de la norme
$$
\Vert(\theta,u,\nabla\Pi)\Vert_{G_T}:= \Vert
\theta\Vert_{L^\infty_T(\dot B^{-1}_{N,\infty})} +\Vert
u\Vert_{L^\infty_T(\dot B^{-1}_{N,\infty})} +\nu\Vert u\Vert_{\widetilde
L^1_T(\dot B^{1}_{N,\infty})} +\Vert\nabla\Pi\Vert_{\widetilde
L^1_T(\dot B^{-1}_{N,\infty})}.
$$
Pour montrer l'unicit\'e dans $G_T,$ nous allons
suivre l'approche utilis\'ee dans \cite{DANCHIN}~:
\'etablir une in\'egalit\'e diff\'erentielle sur
$\Vert(\dt,\du,\nabla\dPi)\Vert_{G_T}$
gr\^ace \`a un argument d'interpolation logarithmique
 puis conclure
\`a l'aide du lemme d'Osgood.

Tout d'abord, en appliquant la proposition \ref{p:transport1} et
en utilisant \eqref{eq:Minkowski}, on obtient
$$
 \Vert\dt\Vert_{L^\infty_t(\dot B^{-1}_{N,\infty})} \leq
Ce^{C\Vert u_2\Vert_{L^1_t(\dot B^{1}_{\infty,1})}}
\int_0^t \Vert\du\,\theta_1\Vert_{\dot B^{0}_{N,\infty}}\,d\tau.
$$
Sachant que $\dot B^0_{N,1}\hookrightarrow
L^N\hookrightarrow\dot B^0_{N,\infty}$ et que $\dot B^0_{\infty,1}
\hookrightarrow L^\infty,$  on peut
\'ecrire~:
$$\begin{array}{lll}
 \Vert\du\,\theta_1\Vert_{\dot B^{0}_{N,\infty}}
&\lesssim&\Vert\du\Vert_{L^\infty}\Vert\theta_1\Vert_{L^N},\\
&\lesssim&\Vert\du\Vert_{\dot B^0_{\infty,1}}
\Vert\theta_1\Vert_{\dot B^0_{N,1}}.
\end{array}$$
Cela permet de conclure que pour tout $t\leq T,$ on a
\begin{equation}\label{densite}
\Vert\dt\Vert_{L^\infty_t(\dot B^{-1}_{N,\infty})} \leq
Ce^{C\Vert u_2\Vert_{L^1_t(\dot B^{1}_{\infty,1})}} \Vert\delta\!
u\Vert_{L^1_t(\dot B^{0}_{\infty,1})} \Vert
\theta_1\Vert_{L^\infty_t(\dot B^{0}_{N,1})}.
\end{equation}
Pour majorer $\du,$ on utilise la remarque \ref{r:stokes} et l'on
obtient pour une constante $C$ ind\'ependante de $t$~:
$$ \dU(t)
\leq C\Bigl(\Vert
\div(u_2\otimes\du)\Vert_{\widetilde L^1_t(\dot B^{-1}_{N,\infty})}
+\Vert\div(\du\otimes u_1)\Vert_{\widetilde
L^1_t(\dot B^{-1}_{N,\infty})} +\Vert\dt\Vert_{\widetilde
L^1_t(\dot B^{-1}_{N,\infty})}\Bigr) $$
avec $\ \displaystyle{\dU(t):=\Vert\du\Vert_{L^\infty_t(\dot B^{-1}_{N,\infty})}
+\Vert\du\Vert_{\widetilde L^1_t(\dot B^{1}_{N,\infty})}
+\Vert\nabla\dPi\Vert_{\widetilde
L^1_t(\dot B^{-1}_{N,\infty})}}.$
\smallbreak
Les termes non lin\'eaires se traitent
\`a l'aide des d\'ecompositions de Bony suivantes~:
$$\begin{array}{lll}
\div(u_2\otimes\du)&=&\dot T_{u^{j}_2}\d_j\du
+\d_j\dot R(u^{j}_2,\du)+\dot T_{\d_j\du}u^{j}_2,\\
\div(\du\otimes u_1)&=&\dot T_{\d_ju_1}\du^j
+\d_j\dot R(u_1,\du^j)
+\d_j\dot T_{\du^j}u_1.
\end{array}
$$
En utilisant les r\'esultats de continuit\'e
de la proposition \ref{p:paraproduit}
(qui se g\'en\'eralisent de fa\c con imm\'ediate
aux espaces $\tilde L_T^\rho(\dot B^s_{p,r})$
voir par exemple \cite{D} ou \cite{DANCHIN}), et \eqref{eq:Minkowski},
on obtient les majorations
$$
\begin{array}{lll}
\Vert\div(u_2\otimes\du)\Vert_{\widetilde
L^1_t(\dot B^{-1}_{N,\infty})} &\!\!\!\lesssim\!\!\! &
\Vert\du\Vert_{\widetilde
L^2_t(\dot B^{0}_{N,\infty})}
\Vert u_1\Vert_{\widetilde L^2_t(\dot B^{0}_{\infty,1})},\\[1ex]
\Vert\div(\du\otimes u_1)\Vert_{\widetilde
L^1_t(\dot B^{-1}_{N,\infty})} &\!\!\!\lesssim\!\!\! &
\Int_0^t\Vert\du\Vert_{\dot B^{-1}_{N,\infty}}
\Vert u_1\Vert_{\dot B^{1}_{\infty,1}}\,d\tau.
\end{array}
$$
En cons\'equence, on a
$$
\dU(t)\lesssim  \Vert u_2\Vert_{\widetilde
L^{2}_t(\dot B^{0}_{\infty,1})}
\Vert\du\Vert_{\tilde L^2_t(\dot B^{0}_{N,\infty})}
+\int_0^t\Bigl(\Vert\dt\Vert_{\dot B^{-1}_{N,\infty}}
+\Vert\du\Vert_{\dot B^{-1}_{N,\infty}}
\Vert u_1\Vert_{\dot B^{1}_{\infty,1}}\Bigr)\,d\tau.
$$
Le premier terme du membre de droite peut \^etre absorb\'e \`a temps
suffisamment petit.
En appliquant l'in\'egalit\'e de Gronwall, on obtient donc pour
$t$ assez petit,
\begin{equation}\label{eq:W1}
\dU(t)\leq C_T\int_0^t\Vert\dt\Vert_{\dot B^{-1}_{N,\infty}}\,d\tau
\end{equation}
pour une constante $C_T$ ne d\'ependant que de normes
de $u_1$ et $u_2$ sur $[0,T].$

 Par ailleurs, d'apr\`es
un r\'esultat d'interpolation logarithmique \'etabli dans la
proposition 1.8 de \cite{D}, on~a
$$
\begin{array}{lll}
\|\du\|_{L^1_t(\dot B^{0}_{\infty,1})} &\lesssim& \|\du\|_{\widetilde L^1_t(\dot B^{0}_{\infty,\infty})} \log
\bigg(e+\frac{\|\du\|_{L^1_t(\dot B^{-1}_{\infty,1})} +\|\du\|_{L^1_t(\dot B^{1}_{\infty,1})}}
 {\|\du\|_{\widetilde L^1_t(\dot B^{0}_{\infty,\infty})}}\bigg),\\[2ex]
 &\lesssim& \|\du\|_{\widetilde L^1_t(\dot B^{1}_{N,\infty})} \log
\bigg(e+\frac{\|\du\|_{L^1_t(\dot B^{1}_{\infty,1})}+t\|\du\|_{L^\infty_t(\dot B^{-1}_{\infty,1})}}
 {\|\du\|_{\widetilde L^1_t(\dot B^{1}_{N,\infty})}}\bigg).
\end{array}
$$
 Par hypoth\`ese, la
fonction  $V$ d\'efinie par
$$
V(t):=\|\du\|_{L^1_t(\dot B^{1}_{\infty,1})}
+t\|\du\|_{L^\infty_t(\dot B^{-1}_{\infty,1})} $$
est born\'ee sur $[0,T].$

En injectant (\ref{densite}) dans \eqref{eq:W1},
et en utilisant l'in\'egalit\'e d'interpolation logarithmique,
on obtient finalement compte tenu de la croissance de la fonction
$x\longmapsto x\ln(e+\frac{y}{x})$ et quitte \`a augmenter $C_T,$
$$ \dU(t) \leq C_T\int_0^t
\dU(s)\log\bigg(e+\frac{V(s)}{\dU(s)}\bigg)\,ds.
$$
Le lemme  d'Osgood  assure   que $\dU(t)=0,$ et donc l'unicit\'e
 sur un petit
intervalle de temps. De l'unicit\'e
locale, on passe \`a  l'unicit\'e sur $[0,T]$ par des arguments
classiques.
\smallbreak Il s'agit maintenant de justifier que
 $(\dt,\du,\nabla\dPi)$ appartient bien \`a l'espace $G_T.$

Vu les hypoth\`eses, la fonction $\dt$ appartient \`a
$L^\infty([0,T];\dot B^{-1}_{N,\infty}).$ En effet, sachant que
$L^{\frac N3,\infty}\hookrightarrow\dot B^{-2}_{N,\infty}$ (cf lemme
\ref{l:lorentz}) et que $\theta_i$ appartient aussi \`a
$L_T^\infty(\dot B^0_{N,1}),$ on d\'eduit
par interpolation que
\begin{equation}\label{eq:u2}
\theta_i\in L_T^\infty(\dot B^{-1}_{N,\infty})\quad\text{pour}\quad i=1,2.
\end{equation}
 Reste \`a justifier  l'appartenance de $(\du,\nabla\dPi)$  \`a
$$
\bigg(\cC\big([0,T];\,\dot B^{-1}_{N,\infty}\big)\cap \widetilde
L^1_T\big(\dot B^{1}_{N,\infty}\big)\bigg)^N \times\bigg(\widetilde
L^1_T\big(\dot B^{-1}_{N,\infty}\big)\biggr)^N.
$$
Pour cela, on  peut d\'ecomposer $u_i$ (que l'on note dor\'enavant
$u$ pour simplifier) en solution libre $u_L:=e^{\nu t\Delta}u_0$
 et fluctuation $\overline u:=u-u_L.$
On voit que $\overline u$ v\'erifie
\begin{equation}\label{eq:unicite}
\left\{\begin{array}{l} \d_t\ov u-\nu\Delta\ov u=
\cP\bigl(\theta\, e_N-\div(u\otimes  u)\bigr),\\
\ov u_{|t=0}=0.
\end{array}\right.
\end{equation}
Tout d'abord, comme $u\in L_T^2(\dot B^0_{\infty,1})$
(par interpolation) et $\dot B^0_{\infty,1}\hookrightarrow L^\infty,$
et comme par ailleurs on a $u\in L_T^\infty(L^{N,\infty}),$
on obtient
$u\in L_T^2(L^{2N}).$
Cela assure que $u\otimes u\in L_T^1(L^N)$ puis, par injection
et diff\'erentiation,
que
\begin{equation}\label{eq:u1}
\div(u\otimes u)\in L_T^1(\dot B^{-1}_{N,\infty}).
\end{equation}

En combinant \eqref{eq:u1}, \eqref{eq:u2}
avec \eqref{eq:Minkowski},
on constate donc que $\ov u$ est solution de l'\'equation
de la chaleur avec donn\'ee initiale nulle
et terme de force dans $\tilde L_T^1(\dot B^{-1}_{N,\infty}).$
La proposition \ref{p:chaleur} assure que
$$
\ov u\in L^\infty_T(\dot B^{-1}_{N,\infty})\cap \tilde
L^1_T(\dot B^{1}_{N,\infty}) \quad\text{et}\quad \nabla\ov\Pi\in\tilde
L_T^1(\dot B^{-1}_{N,\infty}).
$$
Cela ach\`eve la d\'emonstration de la proposition.
\end{p}

La d\'emonstration de l'existence globale
dans le th\'eor\`eme \ref{th:deux}
repose sur des estimations a priori ad\'equates
qui font l'objet de la partie suivante.


\subsubsection{Estimations a priori}\label{s:apriori}

Consid\'erons une solution r\'eguli\`ere $(\theta,u,\nabla\Pi)$
du syst\`eme de Boussinesq. Remarquons tout d'abord
 que la vitesse $u$ v\'erifie
$$u(t)=e^{\nu t\Delta} u_0+\int_0^te^{\nu(t-s)\Delta}\cP
\div(u\otimes u)(s)\,ds+\int_0^t e^{\nu(t-s)\Delta}\cP (\theta
(s)e_N)\,ds.$$ Donc, en vertu des lemmes \ref{l:lorentz2} et
\ref{l:lorentz1}, et sachant que $\Vert u\otimes u\Vert_{L^{\frac
N2,\infty}}\leq \Vert u\Vert_{L^{N,\infty}}^2,$ on~a
$$
\| u\|_{L^\infty_t(L^{N,\infty})} \leq C\Bigl(\|u_0\|_{L^{N,\infty}}
+\nu^{-1}\Vert u\Vert_{L_t^\infty(L^{N,\infty})}^2+
\nu^{-1}\|\theta\|_{L^\infty_t(L^{\frac N3,\infty})}\Bigr).
$$
Comme par ailleurs (cf proposition \ref{p:transportLp}),
$$
\|\theta\|_{L^\infty_t(L^{\frac N3,\infty})}= \|\theta_0\|_{L^{\frac
N3,\infty}},
$$
on en conclut que
\begin{equation}\label{eq:usmall}
\| u\|_{L^\infty_t(L^{N,\infty})}\leq c\nu
\end{equation}
 pourvu que
\begin{equation}\label{eq:smalldata}
\nu^{-1}\|\theta_0\|_{L^{\frac N3,\infty}}
+\|u_0\|_{L^{N,\infty}}\leq c'\nu\end{equation}
avec $c'$ suffisamment petit.
\smallbreak
Il s'agit maintenant de  propager la r\'egularit\'e
n\'ecessaire  pour pouvoir \'etablir l'existence globale
et l'unicit\'e.
Tout d'abord, les lemmes \ref{l:stokes}
et   \ref{l:stokeslimit}, et le fait que
 $L^{N,\infty}\hookrightarrow\dot C^{-1}$ (voir le lemme
 \ref{l:lorentz})  assurent que
 \begin{equation}\label{eq:glob1}
U(t):=\|u\|_{L^\infty_t(\dot B^{-1+\frac Np}_{p,1})}
+\nu\|u\|_{L^1_t(\dot B^{1+\frac Np}_{p,1})}\leq
C\Bigl(\|u_0\|_{\dot B_{p,1}^{-1+\frac Np}}
+\|\theta\|_{L_t^1(\dot B_{p,1}^{-1+\frac Np})}\Bigr)
\end{equation}
pour tout $p\in[N,\infty]$ pourvu que
\eqref{eq:usmall} soit satisfaite.
Pour cela, il suffit que
les donn\'ees initiales v\'erifient \eqref{eq:smalldata},
ce que nous supposerons d\'esormais.
\smallbreak
 D'autre part, en vertu
de  la proposition \ref{p:transport2},
$$\|\theta(t)\|_{\dot B^0_{N,1}}\leq
C\|\theta_0\|_{\dot B^0_{N,1}}
\bigl(1+\|u\|_{L^1_t(\dot B^2_{N,1})}\bigr).$$
Sachant que
 $\dot B^0_{N,1}\hookrightarrow\dot B^{-1+\frac Np}_{p,1},$
on obtient donc
$$U(t)\leq
C\|u_0\|_{\dot B^{-1+\frac Np}_{p,1}}+Ct\|\theta_0\|_{\dot B^0_{N,1}}
+C\nu^{-1}\|\theta_0\|_{\dot B^0_{N,1}}\int_0^tU(\tau)\,d\tau,
$$
puis gr\^ace au lemme de Gronwall,
$$
U(t)\leq C\|u_0\|_{\dot B^{-1+\frac Np}_{p,1}}
e^{Ct\nu^{-1}\|\theta_0\|_{\dot B^0_{N,1}}}
+\nu\Bigl(e^{Ct\nu^{-1}\|\theta_0\|_{\dot B^0_{N,1}}}-1\Bigr).
$$


\subsubsection{D\'emonstration de l'existence globale dans le cas $p=N$}

Soit  $T^*$ le temps d'existence de la solution
$(\theta,u,\nabla\Pi)$  d\'efinie dans le th\'eor\`eme~\ref{th:quatre}.
 Sachant que
$\theta_0\in L^{\frac N3,\infty}$ et que $\theta$ est transport\'e
par un champ de vecteurs qui est dans $L^1_{loc}(0,T^*;{\rm Lip}),$
on~a  $\theta\!\in\!
L^\infty(0,T^*;L^{\frac N3,\infty})$ avec norme constante. Comme de
plus  $u_0\in L^{N,\infty},$ les lemmes \ref{l:lorentz2} et
\ref{l:lorentz1} permettent d'affirmer que $u\in
L^\infty_{loc}(0,T^*;L^{N,\infty}).$

D'autre part, en utilisant le sch\'ema it\'eratif (\ref{schema}), on
obtient une estimation uniforme dans $L^\infty(0,T^*;L^{N,\infty})$.
En effet, on~a pour tout $t\in[0,T^*[,$
$$\|u^{n+1}(t)\|_{L^{N,\infty}}\leq
C\bigl(
\|u_0\|_{L^{N,\infty}}+\nu^{-1}\|\theta_0\|_{L^{\frac N3,\infty}}
+\nu^{-1}\|u^n\|_{L_t^\infty(L^{N,\infty})}^2\bigr)$$
ce qui implique par r\'ecurrence  que
$\|u^n\|_{L_t^\infty(L^{N,\infty})}<2C(\|u_0\|_{L^{N,\infty}}
+\nu^{-1}\|\theta_0\|_{L^{\frac N3,\infty}})<2Cc\nu$ si \eqref{eq:smalldata}
est satisfaite.
Par passage \`a la limite on obtient
$\|u\|_{L^\infty(0,T^*;L^{N,\infty})}\leq2Cc\nu.$

Supposons par l'absurde que $T^*$ est fini. Alors
  les estimations a priori de la section \ref{s:apriori} assurent que
$$\theta\in L^\infty_{T^*}(\dot B^0_{N,1})
\ \text{ et }\ u\in L^\infty_{T^*}(\dot B^{0}_{N,1})\cap
L^1_{T^*}(\dot B^2_{N,1}).$$
Le  th\'eor\`eme \ref{th:quatre} permet alors de prolonger
la solution au-del\`a de $T^*,$ ce qui contredit la maximalit\'e
de  $T^*.$ Cela ach\`eve la d\'emonstration du
th\'eor\`eme \ref{th:deux} dans le cas $p=N.$


\subsubsection{D\'emonstration de l'existence globale dans le cas $p>N$}

La strat\'egie est des plus simples~:
r\'egulariser et tronquer la vitesse initiale afin de se ramener au
cas $p=N$ que l'on sait traiter, puis passer \`a la limite.

Pour la r\'egularisation des donn\'ees, nous ferons appel au lemme
suivant~:
\begin{lem} Soit $\psi\in{\cal C}_c^\infty(\R^N)$
et $(s,p)$ tels que $1\leq p\leq\infty$ et $-N/p'<s\leq N/p.$ Pour
$R>0,$ notons $\psi_R=\psi(\frac\cdot{R}).$ Il existe une
constante $C$ ind\'ependante de $R$ telle que pour toute fonction
$f$ de $\dot B^s_{p,1},$ on ait
$$
\Vert{\psi_Rf}\Vert_{\dot B^s_{p,1}}\leq C\Vert{f}\Vert_{\dot B^s_{p,1}}.
$$
\end{lem}
\begin{p}
Vu les hypoth\`eses sur $s$ et $p,$ la proposition
\ref{p:paraproduit} donne
$$
\Vert{\psi_Rf}\Vert_{\dot B^s_{p,1}}\leq
C\bigl(\Vert\psi_R\Vert_{\dot B^{\frac Np}_{p,1}}
+\Vert{\psi_R}\Vert_{\dot B^{\frac N{p'}}_{p',1}}\bigr) \Vert
f\Vert_{\dot B^s_{p,1}}.
$$
Les normes de $\psi_R$ consid\'er\'ees sont bien finies car ${\cal
C}_c^\infty$ est un sous-espace de $\dot B^{\frac Nr}_{r,1}$  pour tout
$r\in[1,+\infty].$ Par ailleurs, les normes dans ces espaces sont
invariantes par changement d'\'echelle, et ne d\'ependent donc pas
de $R.$
\end{p}
\begin{enumerate}
\item \underline{R\'egularisation des donn\'ees et r\'esolution.}

On pose $u_0^n=\psi(n^{-1}\cdot)\sum_{|q|\leq n}\dot\Delta_qu_0$ o\`u
$\psi$ est une fonction de ${\cal C}_c^\infty(B(0,2))$ \`a valeurs
dans $[0,1]$ et \'egale \`a $1$ sur $B(0,1).$ Il est clair que
$u_0^n$ converge vers $u_0$ au sens des distributions.
 De plus, d'apr\`es le lemme ci-dessus, il
existe une constante $C$ telle que
$$
\forall n\in\N,\;\Vert u_0^n\Vert_{\dot B^{\frac Np-1}_{p,1}} \leq C\Vert
u_0\Vert_{\dot B^{\frac Np-1}_{p,1}},
$$
et l'on v\'erifie facilement en utilisant des in\'egalit\'es de
convolution que
$$
\forall n\in\N,\;\Vert u_0^n\Vert_{L^{N,\infty}} \leq C\Vert
u_0\Vert_{L^{N,\infty}}.
$$
En appliquant le th\'eor\`eme \ref{th:deux} avec $p=N,$
on constate que  si $\nu^{-2}\Vert
\theta_0\Vert_{L^{\frac N3,\infty}} +\nu^{-1}\Vert
u_0\Vert_{L^{N,\infty}}$ est suffisamment petit alors on peut alors
r\'esoudre globalement le syst\`eme de Boussinesq avec donn\'ees
$(\theta_0^n,u_0^n).$
La solution  $(\theta^n,u^n,\nabla\Pi^n)$
ainsi obtenue appartient \`a
$$
{\cal C}(\R_+;\dot B^0_{N,1})\times \Bigl({\cal C}(\R_+;\dot B^0_{N,1})\cap
L^1_{loc}(\R_+;\dot B^2_{N,1}\Bigr)^N
\times\Bigl(L^1_{loc}(\R_+;\dot B^0_{N,1}\Bigr)^N.
$$
\item\underline{Estimations a priori globales.}

Par injection dans les espaces de Besov, $(\theta^n,u^n,\nabla\Pi^n)$
appartient aussi \`a l'espace
$$
E^p:={\cal C}(\R_+;\dot B^0_{N,1})\times \Bigl({\cal C}(\R_+;\dot B^{\frac
Np-1}_{p,1})\cap L^1_{loc}(\R_+;\dot B^{\frac Np+1}_{p,1}\Bigr)^N
\times\Bigl(L^1_{loc}(\R_+;\dot B^{\frac Np-1}_{p,1}\Bigr)^N.
$$
Ainsi les estimations a priori de la partie \ref{s:apriori}
s'appliquent telles quelles. On en d\'eduit des bornes uniformes sur
$(\theta^n,u^n,\nabla\Pi^n)$ dans l'espace $E^p.$
\item \underline{Convergence.}

La convergence de $(\theta^n,u^n,\nabla\Pi^n)_{n\in\N}$ repose sur
des arguments de compacit\'e qui peuvent \^etre d\'eduits des
propri\'et\'es des d\'eriv\'ees temporelles.

Fixons un temps $T>0.$
Par construction, on a  $$
\d_t\theta^n=-\div(u^n\theta^n).$$
Sachant que la suite  $(u^n)_{n\in\N}$ est
 born\'ee dans $L^2_T(\dot B^{\frac
Np}_{p,1})\hookrightarrow L^2_{T}(L^\infty)$ et que la suite
$(\theta^n)_{n\in\N}$ est born\'ee dans
$L^\infty_T(\dot B^0_{N,1})\hookrightarrow L^\infty_T(L^N),$
on en d\'eduit que
$(\d_t\theta^n)_{n\in\N}$ est  born\'ee dans
$L^2_T(W^{-1,N}).$ Cela permet d'affirmer que
$(\theta^n)_{n\in\N}$ est born\'ee dans $\cC([0,T];L^N)$
et \'equicontinue de $[0,T]$
dans $W^{-1,N}.$ Sachant que l'injection
de $L^N$ dans $W^{-1,N}$ est localement compacte,
on peut en d\'eduire en combinant th\'eor\`eme
d'Ascoli et proc\'ed\'e d'extraction diagonal
 la convergence d'une sous-suite de
$(\theta^n)_{n\in\N}$ vers une distribution  $\theta$ qui appartient
localement \`a $\cC([0,T];W^{-1,N}).$

De m\^eme, un examen attentif de l'expression de
$(\d_tu^n)_{n\in\N}$ permet d'affirmer que $(u^n)_{n\in\N}$ est
\'equicontinue de $[0,T]$ dans $B^{-2}_{N,\infty}$
 puis de conclure \`a la convergence \`a
extraction pr\`es. Les d\'etails sont laiss\'es au lecteur.

\item \underline{Conclusion.}

Tout d'abord, les bornes uniformes sur $(\theta^n,u^n)$ permettent
d'affirmer que $(\theta,u)$ appartient en fait \`a
$$
L^\infty_{loc}(\R_+;\dot B^0_{N,1}\cap L^{\frac N3,\infty})
\times\Bigl(L^\infty_{loc}(\R_+;\dot B^{\frac Np-1}_{p,1}\cap
L^{N,\infty}) \cap L^1_{loc}(\R_+;\dot B^{\frac Np+1}_{p,1})\Bigr)^N.
$$
Notons que, m\^eme dans le cas $N=3$ o\`u il faut remplacer
$L^{\frac N3,\infty}$ par $L^1,$ les bornes uniformes dans l'espace de Besov
$\dot B^0_{N,1}$ (qui s'injecte contin\^ument dans $L^N$) emp\^echent
d'\'eventuelles concentrations pour $\theta.$

Par interpolation entre les bornes uniformes
et les r\'esultats de convergence, on peut montrer que
la suite $(\theta^n,u^n)_{n\in\N}$ converge  en un
sens suffisamment fort pour pouvoir passer \`a la limite dans tous
les termes de l'\'equation.

Enfin, sachant que $(\theta,u)$ est solution, il est facile de
montrer que $\d_t\theta+u\cdot\nabla\theta\in L^1_{loc}(\R_+;\dot B^0_{N,1})$
 et que $\d_tu-\nu\Delta u\in L^1_{loc}(\R_+;\dot B^{\frac
Np-1}_{p,1}),$ ce qui assure l'appartenance de $(\theta,u)$
\`a ${\cal
C}(\R_+;\dot B^0_{N,1}) \times{\cal C}(\R_+;\dot B^{\frac Np-1}_{p,1}).$
\hfill\rule{2.1mm}{2.1mm}
\end{enumerate}


\section{Unicit\'e des solutions
d'\'energie finie en dimension deux}

Cette section est consacr\'ee \`a la d\'emonstration
du th\'eor\`eme \ref{th:trois}. 
Rappelons que pour n'importe quelle donn\'ee 
$(\theta_0,u_0)$ dans 
$L^2(\R^2)$ (avec $\div u_0=0$), 
le th\'eor\`eme \ref{solutions-faibles-energie-finie}
assure l'existence globale d'au moins une solutions faible d'\'energie
finie pour le syst\`eme de Boussinesq.
Nous nous proposons de d\'emontrer l'unicit\'e d'une telle solution.
\smallbreak
Dans un premier temps, on va \'etablir que le champ de vitesse
appartient n\'ecessairement  \`a l'espace $\tilde
L^1_{loc}(\R_+;H^2(\R^2)).$ Cette propri\'et\'e 
remarquable (qui nous permettra  d'utiliser
les estimations avec perte   de r\'egularit\'e \'enonc\'ees
dans la proposition \ref{loglip} afin d'\'etablir l'unicit\'e)
a \'et\'e \'etablie par J.-Y. Chemin et N. Lerner dans \cite{CL}
pour le syst\`eme de Navier-Stokes standard. 
L'objet du lemme ci-dessous est de d\'emontrer
que ce r\'esultat demeure pour le syst\`eme de Boussinesq en 
dimension deux.
\begin{lem}\label{l:smoothing}
Soit $(\theta,u,\nabla\Pi)$ une solution de \eqref{eq:boussinesq}
sur $\R_+\times\R^2$ appartenant \`a l'espace d'\'energie~:
\begin{equation}\label{eq:energie}
\theta\in L^\infty(\R_+;L^2(\R^2))\quad\text{et}\quad
u\in L^\infty_{loc}(\R_+;L^2(\R^2))\cap L^2(\R_+;H^1(\R^2)),
\end{equation}
Alors on a $u\in\tilde L^1_{loc}(\R_+;H^2(\R^2))\cap\tilde\cC_{loc}(\R_+;L^2).$
\end{lem}
\begin{dem}
En appliquant l'op\'erateur de projection $\cP$
\`a l'\'equation de la vitesse, on constate que $u$ v\'erifie l'\'equation de
la chaleur~: \begin{equation*}
\begin{cases}
\partial_t u-\nu\Delta u=f_1+f_2\\
u|_{t=0}=u_0\in L^2(\R^2)
\end{cases}\qquad\text{avec}\quad
f_1=-\cP\div(u\otimes u)\quad\text{et}\quad f_2=\cP(\theta e_2).
\end{equation*}
Cela nous am\`ene \`a d\'ecomposer  $u$ en $u_1+u_2$ avec
$$
u_1=\int_0^te^{\nu(t-\tau)\Delta}f_1(\tau)\,d\tau\quad
\text{et}\quad
u_2=e^{t\nu\Delta}u_0+\int_0^te^{\nu(t-\tau)\Delta}f_2(\tau)\,d\tau.
$$
Il est clair que $f_2\in L^1_{loc}(\R_+;L^2).$
De plus, par interpolation, on a visiblement
$u\in L^{\frac83}_{loc}(\R_+;H^{\frac34}).$
En utilisant les lois de produit dans les espaces de
Besov non homog\`enes
(qui sont identiques \`a celles que nous
avons \'enonc\'ees sous la proposition
\ref{p:paraproduit} pour les espaces homog\`enes), on en d\'eduit que
$u\otimes u$ appartient
\`a $L^{\frac43}_{loc}(\R_+;B^{\frac12}_{2,1})$
puis que $f_1\in L^{\frac43}_{loc}(\R_+;B^{-\frac12}_{2,1}).$

Finalement, la proposition \ref{p:chaleur} et la relation \eqref{eq:Minkowski}
garantissent  donc que
\begin{equation}\label{eq:energie1}
u_1\in L^{\frac43}_{loc}(\R_+;B^{\frac32}_{2,1})
\cap\cC(\R_+;B^0_{2,1})\quad\text{et}\quad
u_2\in\tilde L^1_{loc}(\R_+;H^2)
\cap\tilde\cC_{loc}(\R_+;L^2).
\end{equation}
On d\'ecoupe alors $f_1$ en trois  termes~:
$$
f_1=\cP\div(u_1\otimes u_1)+\cP(u_1\cdot\nabla u_2)+
\cP(u_2\cdot\nabla u).
$$
D'apr\`es \eqref{eq:energie1} et par interpolation, on a
$u_1\in L^2_{loc}(\R_+;B^1_{2,1}).$
Sachant que $B^1_{2,1}(\R^2)$ est une alg\`ebre, on en d\'eduit que
$$
\div(u_1\otimes u_1)\in L^1_{loc}(\R_+;B^0_{2,1})
\hookrightarrow L^1_{loc}(\R_+;L^2).
$$
Par ailleurs, en utilisant le fait que $B^1_{2,1}(\R^2)\hookrightarrow
L^\infty(\R^2)$ et que $\nabla u_2\in L^2_{loc}(\R_+;L^2),$
on \'etablit ais\'ement que
$$
u_1\cdot\nabla u_2\in L^1_{loc}(\R_+;L^2).
$$
Enfin, en combinant le fait que
$u_2\in\tilde L^4_{loc}(\R_+;H^{\frac12})$ et que
$\nabla u\in\tilde L^{\frac43}_{loc}(\R_+;H^{\frac12}),$
on trouve que
$$u_2\cdot\nabla u\in L^1_{loc}(\R_+;L^2).$$
Comme l'op\'erateur $\cP$ est continu sur $L^2,$ on peut conclure que
 $f_1\in L^1_{loc}(\R_+;L^2)$
puis, gr\^ace \`a la proposition \ref{p:chaleur}, que  $u_1\in\tilde
L^1_{loc}(\R_+;H^2)\cap\tilde\cC_{loc}(\R_+;L^2).$ \end{dem}
\begin{rem} En \'etant plus soigneux, on peut  \'etablir
que $f_1\in L^1_{loc}(\R_+;B^0_{2,1}),$
ce qui entra\^ine que
$$u_1\in \tilde\cC_{loc}(\R_+;\dot B^0_{2,1})\cap
L^1_{loc}(\R_+;\dot B^{2}_{2,1}).$$
 On retrouve ainsi le fait (bien connu pour
les \'equations de Navier-Stokes) que {\it la fluctuation est plus
r\'eguli\`ere que la solution libre du syst\`eme de Stokes}.
Cependant cette petite am\'elioration ne sera  pas n\'ecessaire pour
d\'emontrer l'unicit\'e dans le th\'eor\`eme \ref{th:trois}.
\end{rem}
Nous sommes maintenant arm\'es  pour d\'emontrer l'unicit\'e.
Consid\'erons donc deux solutions $(\theta_1,u_1,\nabla\Pi_1)$
et $(\theta_2,u_2,\nabla\Pi_2)$ de \eqref{eq:boussinesq}
avec m\^eme donn\'ee initiale
$(\theta_0,u_0)$ dans $L^2(\R^2).$
 Le syst\`eme v\'erifi\'e par la diff\'erence
$(\dt,\du,\nabla\dPi)$ entre ces solutions est
$$
\left\{\begin{array}{l}
\partial_t\dt+\div(u_1\dt)=-\div(\theta_2\,\du),\\
\partial_t\du+\div(u_1\otimes\du)-\nu\Delta\du+\nabla\dPi
=-\div(\du\otimes u_2)+\dt\,e_2.
\end{array}\right.
$$
Le caract\`ere hyperbolique de la
premi\`ere \'equation induit une perte
d'une d\'eriv\'ee dans les estimations de stabilit\'e.
Il va donc falloir d\'emontrer des estimations
au niveau $H^{-1}$ au lieu de $L^2.$
Pour une raison purement technique, et comme l'on dispose d'une grande
 marge de man\oe uvre, il est plus commode
d'op\'erer au niveau $H^{-3/2}$
(mais  n'importe quel exposant strictement compris
entre $-2$ et $-1$ conviendrait aussi).

Pour estimer $\dt$ et $\du,$ nous allons donc appliquer
la proposition \ref{loglip} (ou plut\^ot la remarque \ref{r:loglip}
pour la deuxi\`eme \'equation) avec $s=-3/2$
aux \'equations du syst\`eme.
Pour que cela soit possible, fixons un temps $T>0$ tel que la condition
\eqref{eq:condloglip} de la proposition \ref{loglip}
soit satisfaite par le champ de vitesses $u_1.$
Le lemme \ref{l:smoothing}
 garantit que $u_1\in\tilde L_T^1(H^2)$ et donc l'existence d'un tel  $T.$
En notant
$$\displaylines{
\e_q(t)=C\sum_{q'\leq q}2^{2q'}\Vert\check\Delta_{q'}u_1\Vert_{L_t^1(L^2)}
\qquad
\dTheta(t):=\sup_{\substack{\tau\in[0,t]\\q\geq-1}}
2^{-\frac32q-\e_q(\tau)}
\Vert{\dq\dt(\tau)}\Vert_{L^2},\cr
\dU(t):=\sup_{\substack{\tau\in[0,t]\\q\geq-1}}
2^{-\frac32q-\e_q(\tau)}
\Vert{\dq\du}\Vert_{L^2}+
\nu\sup_{q\geq-1}\int_0^t2^{\frac q2-\e_q(\tau)}
\Vert\dq\du\Vert_{L^2}\,d\tau,}
$$
on obtient alors
$$\begin{array}{lll}
\dTheta(t)&\leq& C\Sup_{q\geq-1}\Int_0^t
2^{-\frac32q-\e_q(\tau)}
\Vert\dq\div(\theta_2\,\du)\Vert_{L^2}\,d\tau,\\[2ex]
\dU(t)&\leq& C(1+\nu t)\Sup_{q\geq-1}\Int_0^t
2^{-\frac32q-\e_q(\tau)}
\bigl(\Vert\dq\div(\du\otimes u_2)\Vert_{L^2}
+\Vert\dq\dt\Vert_{L^2}\bigr)\,d\tau.\end{array}
$$
L'inoffensif  facteur $(1+\nu t)$ provient du fait
que l'on a pris le $\sup$ pour $q\geq-1$ et non pas seulement
pour $q\geq0$ dans le second terme de $\dU(t).$
\smallbreak
Admettons que les termes
de convection se majorent comme suit (voir la
d\'emonstration en appendice)~:
$$
\begin{array}{lll}
\Sup_{q\geq-1}\Int_0^t
2^{-\frac32q-\e_q(\tau)}
\Vert\dq\div\!(\theta_2\,\du)\Vert_{L^2}\,d\tau&\!\!\!\leq\!\!\!& C
\Vert\theta_2\Vert_{L_t^\infty(B^0_{2,\infty})}
\Sup_q\Int_0^t
2^{\frac q2-\e_q(\tau)}
\Vert\dq\du\Vert_{L^2}\,d\tau,\\[1ex]
\Sup_{q\geq-1}\Int_0^t\!
2^{-\frac32q-\e_q(\tau)}
\Vert\dq\div\!(\du\!\otimes\! u_2)\Vert_{L^2}
\,d\tau&\!\!\!\leq\!\!\!& C\Vert u_2\Vert_{\tilde L_t^2(B^1_{2,\infty})}
\Sup_{q\geq-1}\bigl\Vert2^{-\frac q2-\e_q(\tau)}
\Vert\dq\du\Vert_{L^2}\bigr\Vert_{L^2_t}.
\end{array}
$$
En remarquant que pour tout $q\geq-1,$ on~a
$$
\bigl\Vert2^{-\frac q2-\e_q(\tau)}
\Vert\dq\du\Vert_{L^2}\bigr\Vert_{L^2_t}^2
\leq \bigl\Vert2^{\frac q2-\e_q(\tau)}
\Vert\dq\du\Vert_{L^2}\bigr\Vert_{L^1_t}
\bigl\Vert2^{-\frac32q-\e_q(\tau)}
\Vert\dq\du\Vert_{L^2}\bigr\Vert_{L^\infty_t},
$$
on obtient donc
$$\begin{array}{lll}
\dTheta(t)&\leq& C\nu^{-1}
\Vert\theta_2\Vert_{L_t^\infty(B^0_{2,\infty})}
\dU(t),\\[1ex]
\dU(t)&\leq& C(1+\nu t)\Bigl(\nu^{-\frac12}
\Vert u_2\Vert_{\tilde L_t^2(B^1_{2,\infty})}
\dU(t)+\Int_0^t\dTheta(\tau)\,d\tau\Bigr).
\end{array}
$$
Sachant que
$\Vert u_2\Vert_{\tilde L_t^2(B^1_{2,\infty})}
\leq C\Vert u_2\Vert_{L_t^2(H^1)}$
et que
$\Vert\theta_2\Vert_{L_t^\infty(B^0_{2,\infty})}\leq
C\Vert\theta_2\Vert_{L_t^\infty(L^2)}\leq C\Vert\theta_0\Vert_{L^2},$
on en d\'eduit que pour $t$ assez petit, on~a
$$
\dU(t)\leq C(1+\nu t)\int_0^t\dTheta(\tau)\,d\tau\quad\text{et}\quad
\dTheta(t)\leq C\nu^{-1}\Vert\theta_0\Vert_{L^2}\dU(t).
$$
Le lemme de Gronwall permet alors de conclure
\`a l'unicit\'e sur un petit intervalle de temps.

L'unicit\'e globale s'obtient
par un argument de continuit\'e en temps et
de connexit\'e des plus classiques.
\subsubsection*{Fin de la d\'emonstration
du th\'eor\`eme \ref{th:trois}}
Pour montrer que l'\'egalit\'e d'\'energie
\eqref{eq:vitesse} est satisfaite,
il suffit de remarquer que $u$ v\'erifie
le syst\`eme de Navier-Stokes incompressible
en dimension deux avec terme de force dans $L^\infty(\R^+;L^2).$
Enfin, la th\'eorie classique des \'equations
de transport assure
que la norme $L^2$ de $\theta$ est conserv\'ee
au cours de l'\'evolution.
Un argument standard d'analyse fonctionnelle
permet par ailleurs
de montrer que $\theta$ est {\it faiblement} continue
en temps \`a valeurs $L^2.$
Comme la norme $L^2$ est conserv\'ee,
on en d\'eduit que la continuit\'e est forte.
Cela ach\`eve la d\'emonstration du th\'eor\`eme \ref{th:trois}.


\section{Quelques r\'esultats suppl\'ementaires}

\subsection{Un crit\`ere de prolongement}

 On peut d\'emontrer un crit\`ere de
prolongement  bien plus pr\'ecis
que celui qui est \'enonc\'e \`a la fin du th\'eor\`eme
\ref{th:quatre}.
En fait, si 
$$
\int_0^{T^*}\Vert{\nabla u(\tau)}\Vert_{L^\infty}\,d\tau<\infty
$$
alors la solution peut \^etre prolong\'ee au-del\`a 
du temps $T^*.$

Cela r\'esulte du fait que l'int\'egrand apparaissant 
dans  l'in\'egalit\'e  de la proposition 4.7 
peut \^etre remplac\'e par $\Vert\nabla u\Vert_{L^\infty}$
lorsque $|s|<1.$
Il en est donc de m\^eme dans la proposition 4.8. 
En utilisant le th\'eor\`eme 1 de \cite{DAN}, 
on peut alors contr\^oler la norme
$L^\infty(0,T^*;\dot B^0_{N,1})$
de la solution par 
$\int_0^{T^*}\Vert{\nabla u(\tau)}\Vert_{L^\infty}\,d\tau.$
Les d\'etails sont laiss\'es au lecteur.

\subsection{Conditions aux limites p\'eriodiques}

Rappelons qu'il existe une version {\it p\'eriodique}
de la d\'ecomposition de
Littlewood-Paley. De ce fait, tous nos r\'esultats sont
 encore valables {\it mutatis mutandis} dans
le tore $\T^N.$
En fait, sachant que toute solution p\'eriodique \`a moyenne
nulle du syst\`eme de Stokes d\'ecro\^it exponentiellement
en temps, le cas p\'eriodique est m\^eme  plus facile \`a traiter.

\subsection{Termes de forces}

Notons tout d'abord que nos m\'ethodes  n'utilisent jamais la
verticalit\'e du terme $\theta e_N.$
Nos r\'esultats peuvent donc s'adapter au cas plus g\'en\'eral
o\`u $\theta$ est \`a valeurs dans $\R^N.$

 Par ailleurs, pour simplifier la pr\'esentation, nous avons omis les termes
de forces ext\'erieures dans le syst\`eme de Boussinesq.
Indiquons rapidement les hypoth\`eses raisonnables
\`a imposer sur ces termes.

Tout d'abord, si l'on note $\Theta$ le terme de source pour l'\'equation
sur la temp\'erature, et $f$ le terme
de forces ext\'erieures pour l'\'equation sur la vitesse,
on peut obtenir un analogue du th\'eor\`eme
 \ref{solutions-faibles-energie-finie} si l'on suppose en sus que
$$
\Theta\in L^1_{loc}(\R^+;L^p)\quad
\text{et}\quad f\in L^1_{loc}(\R^+;L^2)
\cup L^2_{loc}(\R^+;\dot H^{-1}).
$$
Pour d\'emontrer l'unicit\'e
en dimension deux, il est cependant  imp\'eratif
que $\Theta$ et $f$ soient dans $L^1_{loc}(\R^+;L^2)$
afin que le champ de vitesses appartienne \`a
$\tilde L^1_{loc}(\R^+;H^2).$

\smallbreak
En ce qui concerne les solutions fortes globales,
nos r\'esultats globaux demeurent si l'on suppose que~:
\begin{itemize}
\item  $\Theta$ appartient \`a $L^1_{loc}(\R^+;\dot B^0_{N,1})$
et est petite dans $L^1(\R^+;L^{\frac N3,\infty})$
(resp.  $L^1(\R^+;L^1)$ si $N=3$),
\item $f$ est dans $L^1_{loc}(\R^+;\dot B^{-1+\frac Np}_{p,1})$
et est petite dans $L^p(\R^+;L^{\frac{N(p+2)}{3p},\infty})$ pour un
$p\in[1,\infty]$ (resp.  $L^1(\R^+;L^1)$ si $N=3$ et $p=1$).
\end{itemize}

\subsection{Mod\`ele compl\`etement visqueux}

 Tous nos r\'esultats demeurent valables
sans aucune hypoth\`ese suppl\'ementaire
pour le syst\`eme de Boussinesq compl\`etement visqueux~:
$$
\left\{\begin{array}{l}
\d_t\theta+u\cdot\nabla\theta-\kappa\Delta\theta=0,\\
\d_tu+u\cdot\nabla u-\nu\Delta u+\nabla\Pi=\theta\,e_N,\\
\div u=0,
\end{array}\right.
$$
avec $\kappa$ param\`etre positif
(repr\'esentant la conductivit\'e thermique si $\theta$
est une temp\'erature).

Notre approche permet d'obtenir
des estimations
{\it ind\'ependantes} du param\`etre $\kappa$ sur les solutions construites
et de d\'emontrer ainsi des r\'esultats de convergence pour $\kappa$ tendant
vers $0.$

\subsection{Non existence des singularit\'es
de type jet en dimension trois}

Pour le syst\`eme de Navier-Stokes classique,
il est bien connu que les solutions faibles d'\'energie
finie en dimension trois appartiennent \`a
$L^1_{loc}(\R_+;L^\infty).$

 Pour le  syst\`eme de Boussinesq cette propri\'et\'e
est \'egalement  connue pour des donn\'ees
dans $L^2(\R^3)\times H^1(\R^3)^3,$
ce qui  entra\^ine que 
 les singularit\'es
de type jet (ou {\it squirt} en anglais) ne peuvent pas se produire (voir \cite{CFL}).

Pour les \'equations de Navier-Stokes classiques avec
donn\'ees $L^2,$
on peut  \'etablir une propri\'et\'e
bien plus forte,
\`a savoir que le champ de vitesses
est dans  $L^1_{loc}(\R_+;\dot B^{\frac32}_{2,1})$
ce qui assure notamment l'existence d'un flot continu
en temps et en espace
(voir \cite{Chemin}).
 Cette  propri\'et\'e persiste
pour le syst\`eme de Boussinesq si $3/2<p\leq2$~:
\begin{pro}
En dimension $N=3,$ le champ de vitesses de
toute solution faible d'\'energie finie
du syst\`eme de Boussinesq
avec donn\'ees  $(\theta_0,u_0)\in
L^p(\R^3)\times\bigl(L^2(\R^3)\bigr)^3$
(pour $3/2<p\leq2$)   appartient \`a
$L^1_{loc}(\R_+;B^{\frac32}_{2,1}).$
\end{pro}
\begin{dem}
Nous pr\'esentons une d\'emonstration diff\'erente
de celle de \cite{Chemin}.
Tout d'abord, comme, par construction, le champ de vitesses $u$
 appartient \`a $L^2_{loc}(\R_+;H^1),$
 les inclusions de Sobolev dans les espaces de Lorentz
assurent que $u\in
L^2_{loc}(\R_+;L^{6,2}).$
En utilisant le fait que $\nabla u\in L^2_{loc}(\R_+;L^{2})$
et que $L^2$ co\"{\i}ncide avec l'espace de Lorentz $L^{2,2},$
on~a donc $u\cdot\nabla u\in L^1_{loc}(\R_+;L^{3/2,1}).$
Mais par ailleurs, on a
$$
L^{\frac32,1}=(L^1,L^2)_{\frac23,1}\hookrightarrow
(\dot B^{-\frac32}_{2,\infty},\dot B^0_{2,\infty})_{\frac23,1}
=\dot B^{-\frac12}_{2,1},
$$
donc  $\cP(u\cdot\nabla u)\in L^1_{loc}(\R_+;\dot B^{-1/2}_{2,1}).$
En observant que pour  $3/2<p\leq2,$ on a
 $L^p\hookrightarrow B^{-\frac12}_{2,1},$
on constate que $\cP(\theta\,e_3)\in  L^1_{loc}(\R_+;B^{-1/2}_{2,1}).$
En vertu de   la remarque \ref{r:nonhomogene}, on  conclut alors que
 $$
\int_0^te^{\nu(t-s)\Delta}\cP(\theta\,e_3-u\cdot\nabla u) \,ds \in
L^1_{loc}(\R_+;B^{\frac32}_{2,1}).
$$
Enfin, sachant que $u_0\in L^2,$ on~a
   $e^{t\Delta}u_0\in
L^1_{loc}(\R_+;H^{2-\e})$ pour tout $\e>0.$
Cela ach\`eve la d\'emonstration du r\'esultat.
\end{dem}


\section{Appendice}

Nous d\'emontrons ici
le lemme qui nous a permis d'estimer les termes de convection
dans la preuve de l'unicit\'e du th\'eor\`eme \ref{th:trois}.
\begin{lem}\label{l:quad}
Soit $(\alpha_q)_{q\geq-1}$ une suite
de fonctions positives croissantes sur $[0,T],$
et $\e\in]0,1[.$ On suppose que
pour tout $q'\geq q\geq-1$ et $t\in[0,T],$ on a
\begin{equation}\label{eq:alpha}
0\leq\alpha_{q'}(t)-\alpha_{q}(t)\leq\biggl(\frac{1-\e}2\biggr) (q'-q).
\end{equation}
Alors pour tout  $r\in]1,+\infty],$ il existe une constante
$C$ ne d\'ependant que de $r$ et de $\e$
 telle que l'on ait l'estimation a priori
suivante pour toute fonction $a$ et champ de vecteurs
$b$ \`a divergence nulle sur $\R^2$~:
$$
\sup_{q\geq-1}
\int_0^t2^{-q(1+\e)-\alpha_q(\tau)}
\Vert\dq{\rm div}(ab)\Vert_{L^2}\,d\tau
\leq C\Vert b\Vert_{\tilde L_t^r(B^{\frac2r}_{2,\infty})}
\sup_{q\geq-1}\Bigl\Vert 2^{q(1-\frac2r-\e)-\alpha_q}
\Vert\dq a\Vert_{L^2}\Bigr\Vert_{L_t^{r'}}.
$$
\end{lem}
\begin{dem}
La d\'emonstration repose sur la
d\'ecomposition de Bony, version non homog\`ene.
En tenant compte de $\div b=0,$ on~a (avec la convention
habituelle de sommation sur les indices r\'ep\'et\'es)~:
\begin{equation}\label{eq:bony}
\dq\div(ab)=\dq\bigl(T_{\d_jb}a^j\bigr)
+\dq\bigl(T_{a^j}\d_jb\bigr)
+\dq\d_jR(a^j,b).
\end{equation}
En vertu des propri\'et\'es de localisation
de la d\'ecomposition de Littlewood-Paley, on~a
$$
\dq\bigl(T_{\d_jb}a^j\bigr)=\sum_{|q'-q|\leq4}
\dq\bigl(S_{q'-1}\d_jb\,\Delta_{q'}a^j\bigr).
$$
Ci dessus les op\'erateurs $T$ et $R$ de paraproduit et
de reste  non homog\`enes
se d\'efinissent comme $\dot T$ et $\dot R$
en rempla\c cant chaque  bloc homog\`ene $\ddq$ par
son homologue non homog\`ene $\dq.$

Pour simplifier les calculs,  proc\'edons comme si l'on avait
$\dq\bigl(T_{\d_jb}a^j\bigr)=S_{q-1}\d_jb\,\dq a^j$
(le fait d'avoir \eqref{eq:alpha} permet de justifier
cette approximation).
En utilisant la d\'efinition de $S_{q-1},$
on obtient
$$
\Vert S_{q-1}\d_jb\,\dq a^j\Vert_{L^2}\leq
2^{q(2-\frac2r)}\Vert\dq a\Vert_{L^2}
\sum_{q'\leq q-2}2^{q'(\frac2r-2)}\Vert\Delta_{q'}\nabla b\Vert_{L^\infty}
2^{(q-q')(\frac2r-2)}.
$$
En cons\'equence,  pour tout $0\leq t\leq T,$  on~a
$$\displaylines{
\int_0^t2^{-q(1+\e)-\alpha_q(\tau)}
\Vert S_{q-1}\d_jb\,\dq a^j\Vert_{L^2}\,d\tau
\hfill\cr\hfill\leq\Bigl\Vert 2^{q(1-\frac2r-\e)-\alpha_q}
\Vert\dq a\Vert_{L^2}\Bigr\Vert_{L_t^{r'}}
\sum_{q'\leq q-2}2^{q'(\frac2r-2)}
\Vert\Delta_{q'}\nabla b\Vert_{L_t^r(L^\infty)}
2^{(q-q')(\frac2r-2)}.}
$$
Sachant que $\frac2r-2<0,$ on en d\'eduit alors que
\begin{equation}\label{eq:Tba}\begin{array}[b]{r}
\!\!\!\displaystyle{\sup_{q\geq-1}
\int_0^t2^{-q(1+\e)-\alpha_q(\tau)}
\Vert \dq T_{\d_jb}a^j\Vert_{L^2}\,d\tau}
\\[1.5ex]
\!\!\!\lesssim \Vert
\nabla b\Vert_{\tilde L_t^r(B^{\frac2r-2}_{\infty,\infty})}
\end{array}\sup_{q\geq-1}\Bigl\Vert 2^{q(1-\frac2r-\e)-\alpha_q}
\Vert\dq a\Vert_{L^2}\Bigr\Vert_{L_t^{r'}}.
\end{equation}
De m\^eme, pour estimer
le deuxi\`eme terme de \eqref{eq:bony}, on peut proc\'eder comme
si l'on avait $\dq\bigl(T_{a^j}\d_jb\bigr)=S_{q-1}a^j\,\dq\d_jb.$
 On~a pour tout $0\leq t\leq T,$
$$\displaylines{
2^{-q(1+\e)-\alpha_q(\tau)}
\Vert S_{q-1}a^j\dq\d_jb\Vert_{L^2}\hfill\cr\hfill\leq
2^{q(\frac2r-1)}\Vert\dq\nabla b\Vert_{L^2}
\sum_{q'\leq q-2}2^{-q'(\frac2r+\e)-\alpha_q(\tau)}
\Vert\Delta_{q'}a\Vert_{L^\infty}
2^{(q'-q)(\e+\frac2r)}.}
$$
Sachant que $\alpha_{q}\geq\alpha_{q'}$ pour $q\geq q',$ on en
d\'eduit que
$$\displaylines{
\int_0^t2^{-q(1+\e)-\alpha_q(\tau)}
\Vert S_{q-1}a^j\dq\d_jb\Vert_{L^2}\,d\tau
\hfill\cr\hfill\leq2^{q(\frac2r-1)}\Vert\dq\nabla b\Vert_{L_t^r(L^2)}
\sum_{q'\leq q-2}
\Bigl\Vert 2^{-q'(\frac2r+\e)-\alpha_{q'}}
\Vert\Delta_{q'} a\Vert_{L^\infty}\Bigr\Vert_{L_t^{r'}}
2^{(q'-q)(\e+\frac2r)}.}
$$
Gr\^ace \`a l'in\'egalit\'e de Bernstein, on peut  donc
conclure que
\begin{equation}\label{eq:Tab}\begin{array}[b]{r}
\!\!\!\displaystyle{\sup_{q\geq-1}
\int_0^t2^{-q(1+\e)-\alpha_q(\tau)}
\Vert S_{q-1}a^j\dq\d_jb\Vert_{L^2}\,d\tau\!\!\!}
\\[1.5ex]\lesssim
\Vert \nabla b\Vert_{\tilde L_t^r(B^{\frac2r-1}_{2,\infty})}
\end{array}\sup_{q\geq-1}\Bigl\Vert 2^{q(1-\frac2r-\e)-\alpha_q}
\Vert\dq a\Vert_{L^2}\Bigr\Vert_{L_t^{r'}}.
\end{equation}
Enfin, sachant que
$$
\dq\d_j R(a^j,b)=\sum_{q'\geq q-3}\d_j\dq\bigl(\Delta_{q'}a^j\tilde
\Delta_{q'}b)\quad\text{avec}\quad
\tilde\Delta_{q'}=\Delta_{q'-1}+\Delta_{q'}+\Delta_{q'+1},
$$
et que ${\cal F}\bigl(\d_j\dq\bigl(\Delta_{q'}a^j\tilde
\Delta_{q'}b)\bigr)$ est support\'e dans une boule de taille $2^q,$
on~a en vertu de l'in\'egalit\'e de Bernstein,
$$
2^{-q(1+\e)-\alpha_q(\tau)}\Vert\dq\d_j R(a^j,b)\Vert_{L^2}
\lesssim \sum_{q'\geq q-3}
2^{q(1-\e)-\alpha_q(\tau)}\Vert\Delta_{q'}a\Vert_{L^2}
\Vert\tilde\Delta_{q'}b\Vert_{L^2},
$$
d'o\`u,
$$\displaylines{
\int_0^t2^{-q(1+\e)-\alpha_q(\tau)}\Vert\dq\d_j R(a^j,b)\Vert_{L^2}\,d\tau
\hfill\cr\hfill\lesssim  \sum_{q'\geq q-3}\int_0^t
2^{q'(1-\frac2r-\e)-\alpha_{q'}(\tau)}\Vert\Delta_{q'}a\Vert_{L^2}\,
\bigl(2^{\frac2r{q'}}\!\Vert\tilde\Delta_{q'}b\Vert_{L^2}\bigr)\,
2^{(\alpha_{q'}\!-\!\alpha_q)(\tau)}2^{(q-q')(1-\e)}.}
$$
Pour conclure, l'hypoth\`ese \eqref{eq:alpha} est
cruciale. On obtient alors pour tout $q\geq-1,$
$$
\int_0^t2^{-q(1+\e)-\alpha_q(\tau)}\Vert\dq\d_j R(a^j,b)\Vert_{L^2}\,d\tau
\lesssim
\Vert b\Vert_{\tilde L_t^r(B^{\frac2r}_{2,\infty})}
\sup_{q\geq-1}\Bigl\Vert 2^{q(1-\frac2r-\e)-\alpha_q}
\Vert\dq a\Vert_{L^2}\Bigr\Vert_{L_t^{r'}},
$$
ce qui, joint \`a \eqref{eq:bony}, \eqref{eq:Tba}
et \eqref{eq:Tab}, ach\`eve la d\'emonstration.
\end{dem}

\vskip 9mm
\font\ce=cmtt10
\hspace{-9mm}
\begin{tabular}{lllllllll}
{\sl Rapha\"el Danchin}&&&&&&&&{\sl Marius Paicu}
\\[2ex]
{\sl LAMA, UMR 8050}
&&&&&&&&{\sl Laboratoire de Math\'ematiques}
\\
{\sl
 Universit\'e  Paris-Est,}&&&&&&&&{\sl
 Universit\'e  Paris 11,}
\\
{\sl  61 avenue  du G\'en\'eral de Gaulle,}&&&&&&&&{\sl B\^atiment 425}
\\
 {\sl 94010 Cr\'eteil Cedex, France}&&&&&&&&{\sl 91405 Orsay Cedex, France}
\\
{\ce danchin@univ-paris12.fr}&&&&&&&&
{\ce marius.paicu@math.u-psud.fr}
\end{tabular}

\end{document}